\documentclass[a4paper, leqno]{amsart}

\usepackage[latin1]{inputenc}
\usepackage[english]{babel}
\usepackage{amssymb,exscale,xy,mathrsfs,bbm}\xyoption{all}

\theoremstyle{plain}
\newtheorem{theorem}{Theorem}[section]         
\newtheorem{corollary}[theorem]{Corollary}     
\newtheorem{lemma}[theorem]{Lemma}             
\newtheorem{proposition}[theorem]{Proposition} 

\theoremstyle{definition}
\newtheorem{observation}[theorem]{Key Observation}

\newtheorem{remark}[theorem]{Remark}

\renewcommand{\labelenumi} {(\alph{enumi})}    

\renewcommand{\theenumi} {(\alph{enumi})}

   \newcommand{\norm}[1] {\| #1 \|}
   
   \newcommand{\bignorm}[1]{\bigl\| #1 \bigr\|}
   \newcommand{\Bignorm}[1]{\Bigl\| #1 \Bigr\|}
   
   \newcommand{\biggnorm}[1]{\biggl\| #1 \biggr\|}
   \newcommand {\sfrac}[2] { {{}^{#1}\!\!/\!{}_{#2}}} 
   \newcommand {\einhalb} {\sfrac{1}{2}} 
   \newcommand {\pihalbe} {\sfrac{\pi\,}{2}}
   \renewcommand {\AA} {\mathbb A}
   \newcommand {\RR} {\mathbb R}
   \newcommand {\NN} {\mathbb N}
   \newcommand {\XX} {\mathbb X}
   \newcommand {\ZZ} {\mathbb Z}
   \newcommand {\VV} {\mathbb V}
   \newcommand {\BB} {\mathbb B}
   \newcommand {\PP} {\mathbb P}
   \newcommand {\CC} {\mathbb C}
   \renewcommand {\VV} {\mathbb V}
   \newcommand {\cF} {\mathcal F}
   \newcommand {\RANGE} {\mathcal R}
   \newcommand {\DOMAIN} {\mathcal D}
   \newcommand {\cN} {\mathcal N}       
   \newcommand {\sD} {\mathscr D}
   \newcommand {\Morrey} {\mathcal M}
   \newcommand {\sF} {\mathscr F}
   \newcommand {\eins} {\mathbbm 1}
   \newcommand {\la}{\lambda}
   \newcommand {\si}{\sigma}

   \newcommand {\ga}{\gamma}
   \newcommand {\om}{\omega}
   \newcommand {\Om}{\Omega}
   \newcommand {\al}{\alpha}
   \newcommand {\eps}{\epsilon}
   \newcommand {\Sec}[1] {S({#1})}
   
   \newcommand {\idual}[3] {\langle #1, #2 \rangle_{#3}}
   
   \newcommand {\dual}[2] {\idual{#1}{#2}{} }
   \newcommand {\emb} {\hookrightarrow}
   
   \newcommand {\ess} {{\rm ess.}}
   \newcommand {\suchthat}{:\;}

\begin{document}
  \date{21.12.2007}

  \title[On Kato's method for Navier--Stokes Equations] 
        {On Kato's method for Navier--Stokes Equations} 

  \author[Bernhard H. Haak]{Bernhard H. Haak}
  \address{Institut de Math\'ematiques de Bordeaux\\%
           351 cours de la  Lib\'eration\\33405 Talence CEDEX\\France}
  \email{Bernhard.Haak@math.u-bordeaux1.fr}

  \author[Peer Chr. Kunstmann]{Peer Chr. Kunstmann} 
  \address{Institut f\"ur Analysis\\ Universit\"at %
           Karlsruhe\\Englerstra\ss{}e 2\\ 76128 Karlsruhe\\ Germany}       
  \email{Peer.Kunstmann@math.uni-karlsruhe.de}

  \thanks{This research was done at Karlsruhe University of
    Technology. The authors kindly
    acknowledge support from Deutsche Forschungsgemeinschaft, 
    contract number WE 2847/1-2}

   \subjclass{35Q30, 47D06, 46B70}

   \keywords{Mild solutions, Navier-Stokes equations, Kato's method, %
    parabolic equations, quadratic non-linearity, admissibility of unbounded operators}

\begin{abstract}
We investigate Kato's method for parabolic equations with a 
quadratic non-linearity in an abstract form. We extract several 
properties known from linear systems theory which turn out to be
the essential ingredients for the method. We give necessary and 
sufficient conditions for these conditions and provide new and 
more general proofs, based on real interpolation. 
In application to the Navier-Stokes equations, our approach unifies 
several results known in the literature, partly with different proofs.
Moreover, we establish new existence and uniqueness results for 
rough initial data on arbitrary domains in $\RR^3$ and irregular 
domains in $\RR^n$.
\end{abstract}

  \maketitle

\section{Introduction}

Let $\Om\subset\RR^n$ be a domain, i.e. an open and connected subset.
In this paper we study the Navier--Stokes equation in the form
\begin{equation}\label{eq:nse}\tag{NSE}
 \left.\begin{array}{rcl}
 u_t-\Delta u+(u\cdot\nabla)u+\nabla p &=& f,\quad(t>0)\\ 
              \nabla\cdot u &=& 0 \\
            u(0,\cdot) &=& v_0 \\
           u|_{\partial\Om} &=& 0.
 \end{array}\right\}
\end{equation}
The equation \eqref{eq:nse} describes the motion of an incompressible
fluid filling the region $\Om$ under ``no slip'' boundary conditions, 
where $u=u(t,x)\in\RR^n$ denotes the unknown velocity vector at time
$t$ and point $x$, $p=p(t,x)\in\RR$ denotes the unknown pressure, and
$v_0$ denotes the initial velocity field which is also assumed to be
divergence--free, i.e. $\nabla\cdot v_0=0$. Of course, the boundary
condition is not present in case $\Om=\RR^n$. Observe already that
$\nabla\cdot u=0$ allows to rewrite
$(u\cdot\nabla)u=\nabla\cdot(u\otimes u)$.

Initiated perhaps by {\sc Cannone}'s work (\cite{Cannone:Buch})
there has been  a lot of interest in the last decade in mild solutions of 
\eqref{eq:nse} (see e.g. \cite{Amann:strong-solvability-NS,KochTataru,
KozonoYamazaki,Lemarie:NSE,Meyer:Navier-Stokes,Sawada:time-local}) for
initial data in so-called critical spaces. All these results rely on
variations of {\sc Kato}'s method (\cite{FujitaKato}) which allows
to obtain global 
solutions if the initial data are small by a fixed point argument 
(which is based on Banach's fixed point principle or, equivalently, 
on a direct fixed point iteration).

The fixed point equation is obtained from \eqref{eq:nse} by first applying 
the Helm\-holtz projection $\PP$ to get rid of the pressure term
\begin{equation}\label{eq:nseP} 
 \left.\begin{array}{rcl}
 u_t-\PP\Delta u+\PP\nabla\cdot(u\otimes u) &=& \PP f,\quad(t>0)\\ 
              \nabla\cdot u &=& 0 \\
            u(0,\cdot) &=& v_0 \\
           u|_{\partial\Om} &=& 0.
 \end{array}\right\}
\end{equation}
The operator $-\PP\Delta$ with Dirichlet boundary conditions is, 
basically, the \emph{Stokes operator} $A$ which -- hopefully -- is the negative
generator of a bounded analytic semigroup $T(\cdot)$, the \emph{Stokes semigroup}, 
in the divergence--free function space $X$ under consideration. Then 
the solution to \eqref{eq:nseP} is formally given by the 
variation-of-constants formula 
\begin{equation}\label{eq:fix-point} 
 u=T(\cdot)v_0 - T(\cdot)*\PP\nabla\cdot(u\otimes u) + T(\cdot)*\PP f.
\end{equation}
If one can give sense to the Helmholtz projection $\PP$, the Stokes operator $A$ 
and the Stokes semigroup $T(\cdot)$, this is a fixed point equation for $u$. 
A \emph{mild solution} to \eqref{eq:nse} is a solution to \eqref{eq:fix-point}. 

The non-linearity is quadratic and may be rewritten using the bilinear map 
$F(u,v):=\PP\nabla\cdot(u \otimes v)$. The natural space for a fixed 
point argument yielding global solutions would be $C([0,\infty),X)$, but this 
rarely works for critical spaces. The idea of Kato's method for the critical 
space $X=L^3$ on $\Om=\RR^3$ is to use an auxiliary space $Z=L^q$ with $q\in(3,6]$ 
and a weighted sup-norm with a polynomial weight $t^\al$ and to carry out the 
iteration scheme in a suitable function space with norm
\[
  \norm{t\mapsto u(t)}_{L^\infty(\RR_+,L^3)}
 +\norm{t\mapsto t^{\frac{1}{2}-\frac{3}{2q}}u(t)}_{L^\infty(\RR_+,L^q)}.
\]
In our paper we the note $L^p_\al((0,\tau),X)$ the space of all $X$-valued
measurable functions $f$ such that 
\[
  \norm{f}_{L^p_\al((0,\tau),X)}:=\norm{t\mapsto t^\al f(t)}_{L^p((0,\tau),X)}<\infty.
\]
As {\sc Cannone} observed (\cite{Cannone:Buch}, see also
\cite{KozonoYamazaki}), Kato's approach leads to Besov
spaces in a natural way. On suitable domains $\Om\neq\RR^n$, 
{\sc Amann}'s work (\cite{Amann:strong-solvability-NS}) underlined the
fundamental role of real 
interpolation and of abstract extrapolation and interpolation scales. 
The present paper takes up this point of view.

We start our main results with an abstract version of Kato's method for
parabolic equations with quadratic non-linearity 
(Theorem~\ref{thm:kato-abstrakt}), which clearly isolates the properties 
one has to check for in order to obtain local solutions for arbitrary
data or global solutions for small initial data. These properties
\ref{item:C-zul}, \ref{item:B-zul}, and  
\ref{item:faltung-bdd} only concern linear problems. 

In the literature, there is an abstract version of Kato's method due to 
{\sc Weiss\-ler} \cite{Weissler:NSE}, formulated for parabolic equations with quadratic
non-linearity. The approach, however, is different already for the bilinear term
(see Remark~\ref{rem:weissler}), 
and in extension to Weissler's result we do not only consider weighted sup-norms 
for functions with values in an auxiliary space, but also weighted
$L^p$--spaces with polynomial weights $t^\al$ for $p\in[2,\infty]$ 
(the restriction  $p\ge2$ is due to the quadratic nature of the 
non-linearity). Moreover, in our second main result 
(Theorem~\ref{thm:discussion-A1-A3}) we give necessary 
and sufficient conditions for the properties \ref{item:C-zul}, 
\ref{item:B-zul}, and \ref{item:faltung-bdd}.
We were led to these results by our previous work on linear systems 
of the form
\begin{equation}\label{eq:lin-sys}
 \left.\begin{array}{rcl}
 x'(t)+Ax(t) &=& Bu(t),\quad t>0\\ 
        y(t) &=& Cx(t),\quad t>0 \\
        x(0) & = & x_0
 \end{array}\right\}
\end{equation}
Theorem~\ref{thm:kato-abstrakt} is actually a result on a quadratic 
feedback law $u(t)=F(y(t),y(t))$ for \eqref{eq:lin-sys}. 
In \eqref{eq:lin-sys}, $C$ and $B$ are unbounded linear operators 
(in the application to \eqref{eq:nse} they are the identity 
on suitable spaces, see below), and \ref{item:C-zul} and 
\ref{item:B-zul} simply mean that they are \emph{admissible} in the
sense of linear systems theory for the corresponding weighted Bochner spaces.
The conditions in Theorem~\ref{thm:discussion-A1-A3}~\ref{item:discussion-A1}
and \ref{item:discussion-A2} are generalisations of our results in 
\cite{HaakKunstmann:weighted-Lp-admiss} to the case of not necessarily 
densely defined operators $A$.  Moreover, we give here new and very 
transparent proofs based on real interpolation 
(see Section~\ref{sec:proof-of-kato-discussion}) whereas the proofs in 
\cite{HaakKunstmann:weighted-Lp-admiss} relied on $H^\infty_0$--functional
calculus arguments. 

In Section~\ref{sec:applications} we apply our abstract results to obtain mild solutions 
to \eqref{eq:nse}. On $\RR^n$ we reobtain {\sc Cannone}'s result (\cite{Cannone:Buch})
on initial values in Besov spaces (see Subsection \ref{sec:lebesgue}). 
In Subsection \ref{sec:lorentz} we show that, in close analogy to 
Subsection~\ref{sec:lebesgue}, one may likewise use weak Lebesgue spaces as 
auxiliary spaces $Z$ which leads to mild solutions for initial values in Besov type spaces 
that are based on weak Lebesgue spaces. This result is new.
In Subsection~\ref{sec:morrey} we use Morrey spaces as auxiliary spaces and obtain
results similar to those in {\sc Kozono} and {\sc Yamazaki}
\cite[Theorem 3]{KozonoYamazaki} on initial values in Besov type
spaces based on Morrey spaces. Our approach allows to reproduce their
result even under weaker conditions for the initial value.
In Subsection \ref{sec:hoelder} we give a variant of a result due to
{\sc Sawada} (\cite{Sawada:time-local}) on time-local solutions for initial values in Besov spaces 
$B^{-1+\eps}_{\infty,p}$ with $p\in(n,\infty)$, but with a different proof.
For the quadratic term we simply use the 
product inequality for H\"older continuous functions whereas the keystone of the proof 
in \cite{Sawada:time-local} was a H\"older type inequality for functions in general Besov spaces.
Subsection \ref{sec:domains-L2} studies mild solutions for arbitrary domains 
$\Om\subset\RR^3$, and we improve results due to 
{\sc Sohr} (\cite[Theorem V.4.2.2]{Sohr:Buch}) and 
{\sc Monniaux} (\cite[Theorem 3.5]{Monniaux:arbitrary}). Moreover, our approach allows
to compare both results.
In Subsection \ref{sec:domains-Lq} we assume that Helmholtz projection and Stokes 
semigroup act in a scale of $L^q$-spaces, $q\in[q_0',q_0]$, and investigate how the value
of $q_0>2$ affects existence of mild solutions for certain initial values. 
It turns out that, already under these relatively weak assumptions, a larger
$q_0$ allows for more initial values, where the case $q_0>\max(4,n)$ needs an additional 
gradient estimate for the Stokes semigroup. In any case, these new results make very clear
which properties one has to check for the Stokes semigroup in order to obtain mild solutions
for ``rough'' initial values, i.e. for initial values in suitable extrapolation spaces.
We mention that there are other approaches to the Navier-Stokes
equations for rough initial data or on general domains (see e.g.,
\cite{KatoPonce,KozonoYamazaki,KochTataru,FarwigKozonoSohr}), and we
shall comment on them at the end of each subsection in
Section~\ref{sec:applications}. 

The paper is organised as follows. In Section~\ref{sec:prelim} we collect basic facts on the 
Helmholtz decomposition and the Stokes semigroup for arbitrary domains. Those are the 
basis for applications of the abstract results to \eqref{eq:nse} 
in Section~\ref{sec:applications}. In Section~\ref{sec:abstract-kato} we present our abstract results,
a part of the proofs is relegated to Section~\ref{sec:proof-of-kato-discussion}.
In an appendix we have gathered facts on Besov spaces based on weak
Lebesgue spaces that are needed in Subsection \ref{sec:lorentz} and
facts on Morrey spaces that are needed in Subsection \ref{sec:morrey}.

\medskip

{\bf Acknowledgement:} The authors thank the unknown referee for several
suggestions which helped to 
improve 
the article, in particular for
drawing our attention to \cite{KozonoYamazaki} and interpolation spaces
of Morrey spaces.

\section{Preliminaries}\label{sec:prelim}

Let $n\geq 2$ and let $\Om\subseteq\RR^n$ be an arbitrary open and connected 
subset. We start with basics on the Helmholtz decomposition in $L^q(\Om)^n$
where $q\in(1,\infty)$. To this end we define
\[
 \dot{W}^1_q(\Om) := \bigl\{ [u]= u+\CC \suchthat u\in L^q_{loc}(\Om) 
    \text{ and }\nabla u\in L^q(\Om)^n \bigr\}
\]
with norm $\norm{ u }_{\dot{W}^1_q(\Om)}:=\norm{ \nabla u }_{L^q(\Om)^n}$.
Although the Navier--Stokes equations involve real valued functions, we 
consider complex function spaces here, since our abstract arguments below
shall deal with complex Banach spaces.

The space $\dot{W}^1_q(\Om)$ is a Banach space and the linear map
$\nabla_q:\dot{W}^1_q(\Om)\to L^q(\Om)^n$, $u\mapsto \nabla u$,
is isometric. We also define
\[
 (\dot{W}^1_q(\Om))'
 :=\{\phi \suchthat \dot{W}^1_q(\Om) \to \CC \suchthat \phi\ \mbox{is linear and continuous}\}
\]
with the usual operator norm. Then $(\dot{W}^1_q(\Om))'$ is a Banach space and 
the dual map
$(\nabla_{q'})':L^q(\Om)^n\to(\dot{W}^1_{q'}(\Om))'$ of 
$\nabla_{q'}$ 
is surjective with norm $\leq1$.

We recall the space $\sD(\Om)=C^\infty_c(\Om)$ of test functions
and the dual space $\sD'(\Om)$ of distributions on $\Om$.
\begin{remark}\label{rem:necas}
If $u\in\sD'(\Om)$ satisfies $\nabla u\in L^q(\Om)^n$ then 
$u$ belongs to $L^q_{loc}(\Om)$ (\cite{Necas:methodes-directes}).
\end{remark}

Now let $G^q(\Om):=\mbox{Im}\nabla_q=\nabla_q\dot{W}^1_q(\Om)$ 
denote the space of gradients in $L^q(\Om)^n$ and 
$L^q_\sigma(\Om):=\mbox{Ker}(\nabla_{q'})'$ denote the space of 
divergence-free vector fields.
\begin{remark}\label{rem:orthogonality}
It is clear from the construction that
\begin{eqnarray*}
 G^q(\Om) &=& \{f\in L^q(\Om)^n \suchthat \forall g\in L^{q'}_\si(\Om)
 \suchthat \dual{ f}{g } = 0\}  \quad\mbox{and}\\
 L^q_\si(\Om) &=& \{f\in L^q(\Om)^n \suchthat \forall g\in G^{q'}(\Om)
 \suchthat \dual{ f}{g}=0\}.
\end{eqnarray*}
\end{remark}

Let $\sD_\sigma(\Om) := \{\phi\in\sD(\Om)^n\suchthat\nabla\cdot\phi=0\}$
denote the space of divergence-free test functions.
The following theorem holds.
\begin{theorem}[de Rham]\label{thm:deRham}
Let $T\in\sD'(\Om)^n$. There is an $S\in\sD'(\Om)$
with $T=\nabla S$ if and only if $T$ vanishes on $\sD_\sigma(\Om)$.
\end{theorem}

\begin{remark}
This was first noticed by {\sc Lions} \cite[p.67]{Lions:quelques-methodes}
who resorted to a result due to {\sc de Rham}
\cite[Theorem 17', p.114]{DeRham:varietes-differentiables}.
We refer to \cite{Simon:primitives} for more details and an elementary proof.
\end{remark}

Now we are able to prove the following representation of the space 
$L^q_\sigma(\Om)$ which is often taken as the definition. The 
argument in the proof is the same as in \cite[p.67]{Lions:quelques-methodes}.
\begin{proposition}
For any $q\in(1,\infty)$, the space $L^q_\sigma(\Om)$ is the closure of
$\sD_\sigma(\Om)$ in $L^q(\Om)^n$.
\end{proposition}
\begin{proof}
Let $\phi\in\sD_\sigma(\Om)$ and $u\in\dot{W}^1_{q'}(\Om)$. Then
\[
  \dual{(\nabla_{q'})'\phi}{u}  =\dual{ \phi }{\nabla_{q'} u}  =-\dual{\nabla\cdot\phi}{u}=0,
\]
and $\phi\in L^q_\sigma(\Om)$. To show density of $\sD_\si(\Om)$ in
$L^q_\si(\Om)$ we take $g\in L^{q'}(\Om)^n$ such that $\dual{ g}{\cdot}$
vanishes on $\sD_\si(\Om)$ and have to show that 
$\dual{g}{ \cdot}$ vanishes on $L^q_\si(\Om)$. 
By Theorem~\ref{thm:deRham} we find $v\in\sD'(\Om)$ such that
$\nabla v=g$. By Remark~\ref{rem:necas} we have $v\in L^{q'}_{loc}(\Om)$,
i.e. $g\in G^{q'}(\Om)$. Now we use Remark~\ref{rem:orthogonality}.
\end{proof}

Concerning the Helmholtz projection we quote the following theorem, which is 
the essence of the approach in \cite{SimaderSohr:new-approach}.
\begin{theorem}\label{thm:sim-sohr}
Let $q\in(1,\infty)$. Then $L^q(\Om)^n=L^q_\si(\Om)\oplus G^q(\Om)$
if and only if the operator 
$N_q:=(\nabla_{q'})'\nabla_q:\dot{W}^1_q(\Om) \to (\dot{W}^1_{q'}(\Om))'$
is bijective.
\[
\xymatrix{
    & **[r] (\dot W_{q'}^1(\Om))' \\                  
L^q(\Om)^n \ar[ur]^{(\nabla_{q'})'} &            \\
    & \dot W_q^1(\Om) \ar[ul]_{\nabla_q} \ar[uu]_{N_q}
}
\]
If the operator $N_q$ has a bounded inverse 
$N_q^{-1}:(\dot{W}^1_{q'}(\Om))' \to \dot{W}^1_{q}(\Om)$,
then the projection from $L^q(\Om)^n$ onto $L^q_\si(\Om)$ with
kernel $G^q(\Om)$ is given by $\PP_q:=I-\nabla_q N_q^{-1} (\nabla_{q'})'$.
This projection is called the Helmholtz projection in $L^q(\Om)^n$.
\end{theorem}
\begin{proof}
If $N_q$ is bijective then its inverse $N_q^{-1}$ is bounded by the open 
mapping theorem, the projection $\PP_q$ has the desired properties and 
we obtain $L^q(\Om)^n=L^q_\si(\Om)\oplus G^q(\Om)$.

Conversely, if $L^q(\Om)^n=L^q_\si(\Om)\oplus G^q(\Om)$ then 
$\nabla_q\dot{W}^1_q(\Om)\cap\text{Ker}(\nabla_{q'})'=\{0\}$ and
$N_q$ is injective. 
Moreover $(\dot{W}^1_{q'}(\Om))'=(\nabla_{q'})'G^q(\Om)$, thus 
$N_q$ is surjective by $G^q(\Om)=\nabla_q\dot{W}^1_q(\Om)$.
\end{proof}

In \cite{SimaderSohr:new-approach}, the operator $-N_q$ is interpreted as a 
weak version of the Neumann--Laplacian on $\Om$.
For $q=2$, $N_2$ is always bijective, and $\PP_2$ is the orthogonal
projection from $(L^2(\Om))^n$ onto $L^2_\sigma(\Om)$. This follows from
Remark~\ref{rem:orthogonality} or from Theorem~\ref{thm:sim-sohr} via Lax--Milgram.
\begin{remark}\label{rem:div-null}
Since $\sD(\Om)\subset\dot{W}^1_q(\Om)$, Remark~\ref{rem:orthogonality} shows that
$u\in L^q_\si(\Om)$ implies $\nabla\cdot u=0$ in the sense of distributions.
However $L^q_\si(\Om)$ also contains information on the behaviour of
$u$ at the boundary (see \cite{SimaderSohr:new-approach}, 
\cite[Lemma II.2.5.3]{Sohr:Buch}): for example for bounded Lipschitz 
domains one has $u \in L^q_\si(\Om)$ if and only if $\nabla\cdot u = 0$ on 
$\Om$ and $\nu\cdot u=0$ on $\partial \Om$ 
where $\nu$ denotes the outer normal unit vector.
\end{remark}

\subsection*{Function spaces.}
For $q\in(1,\infty)$, we use the usual notation and write 
$W^1_{q,0}(\Om) = %
\overline{ {\mathscr D}(\Om) }^{\norm{\cdot}_{W^1_q}}$
and $W^{-1}_q(\Om):=(W^1_{q',0}(\Om))'$. Moreover, we let
\[
 \dot{W}^1_{q,0}(\Om):=(W^1_{q,0}(\Om),\norm{ \nabla\cdot }_q)^\sim
 \qquad\mbox{and}\qquad
  \dot{W}^{-1}_q(\Om):=( \dot{W}^1_{q',0}(\Om) )',
\]
where ${}^\sim$ denotes the completion. Then $\dot{W}^{-1}_q(\Om)$ consists 
of all $\phi\in W^{-1}_q(\Om)$ satisfying
\[
 \norm{ \phi }_{\dot{W}^{-1}_q}
 =\sup\{ |\phi(v)| \suchthat v\in W^1_{q',0}(\Om), \norm{\nabla v }_{q'}\le 1\}<\infty.
\]
The corresponding spaces of ``divergence-free'' vectors are
\[
 \VV_q(\Om):=W^1_{q,0}(\Om)^n\cap L^q_\si(\Om)\quad\mbox{and}\quad
 \dot{\VV}_q(\Om):=(\VV_q(\Om),\norm{ \nabla\cdot }_q)^\sim.
\]
\[
 \dot{\VV}^{-1}_q(\Om):=(\dot{\VV}_{q'}(\Om))'
 =\{\phi \suchthat \VV_{q'}(\Om)\to\CC\text{ linear}\!\suchthat
 \phi\,\text{ is continuous for $\norm{ \nabla\cdot }_{q'}$}\,\}
\]
with the natural operator norm. 
Then $\dot{\VV}_q(\Om)$ is a Banach space and $\VV_q(\Om)$ is a Banach space for
the norm of $W^1_q(\Om)^n$ and a dense subset of $\dot{\VV}_q(\Om)$. 
\begin{lemma}
The set $\sD_\si(\Om)$ is dense in $(\VV_q(\Om),\norm{ \cdot }_{W^1_q})$
and in 
$(\dot \VV_q(\Om),\norm{ \nabla\cdot }_q)$.
\end{lemma}
\begin{proof}
By definition it suffices to consider $\VV_q(\Om)$. It is clear that 
$\sD_\si(\Om)\subset \VV_q(\Om)$. Now take $\phi\in(W^1_{q,0}(\Om)^n)'$
such that $\phi$ vanishes on $\sD_\si(\Om)$. Notice that $\phi$ is a 
distribution on $\Om$. By Theorem~\ref{thm:deRham} there exists 
$h\in\sD '(\Om)$ satisfying $\phi=\nabla h$ and $h$ is unique up to a 
constant. Since $\phi$ can be represented as a sum of partial derivatives
of $L^{q'}$-functions, we conclude that we can assume $h\in L^{q'}(\Om)$. 

For $u\in\VV_q$ we choose a sequence $(u_k)$ in $\sD(\Om)^n$ such that
$u_k\to u$ in $W^1_{q,0}(\Om)^n$, and we finally obtain
\[
 \phi(u)=\dual{\nabla h}{u}=\lim_k \dual{\nabla h}{u_k}=
- \lim_k \dual{h}{\nabla\cdot u_k}= -\dual{ h}{\nabla\cdot u}=0
\]
by $\nabla\cdot u=0$ (see Remark~\ref{rem:div-null}). This ends the proof.
\end{proof}

Coming back to the Navier-Stokes equation we notice that,
for $u\in L^q(\Om)^n$, we have $u\otimes u\in L^{\sfrac{q}2}(\Om)^{n\times n}$ 
and $\nabla\cdot(u\otimes u)\in \dot{W}^{-1}_{\sfrac{q}2}(\Om)^n$. 
Applying the Helmholtz projection to get rid of the pressure term 
$\nabla p$ in \eqref{eq:nse} thus needs extensions $P_q$ of the Helmholtz 
projection $\PP_q$ to $\dot{W}^{-1}_q(\Om)^n$, $q\in(1,\infty)$.
Those are defined by restriction (as in, e.g., \cite{Sohr:Buch}, 
\cite{Monniaux:arbitrary}):
\[                                                      
 P_q:\dot{W}^{-1}_q(\Om)^n\to\dot{\VV}^{-1}_q(\Om), \qquad
 P_q\phi(v):=\phi|_{\dot{\VV}_{q'}(\Om)}.
\]
Observe that this is meaningful since 
$\dot{\VV}_{q'}(\Om)\subset\dot{W}^1_{q',0}(\Om)^n$.
Moreover, $P_q$ is linear and continuous. We show that $P_q$ and $\PP_q$
are consistent.
\begin{lemma}
We have $P_q\phi=\PP_q f$ for each $\phi\in\dot{W}^{-1}_q(\Om)^n$ and
$f\in L^q(\Om)^n$ such that $\phi(v)=\dual{ f}{v}$ for all
$v\in W^1_{q',0}(\Om)^n$.
\end{lemma}
\begin{proof}
It suffices to check equality on 
$\VV_{q'}(\Om)=W^1_{q',0}(\Om)^n\cap L^{q'}_\si(\Om)$.
For $v\in\VV_{q'}(\Om)$ we have
\[
 P_q\phi(v)=\phi(v)=\dual{ f}{v}=\dual{ f}{\PP_{q'}v}
=\dual{ \PP_q f}{v}
\]
by $\PP_{q'}v=v$ and $(\PP_{q'})'=\PP_q$.
\end{proof}

\subsection*{The Stokes operator.}
We define the Stokes operator in $L^2_\si(\Om)$ by the form method.
To this end we let $\VV:=\VV_2=L^2_\si(\Om)\cap(W^1_{2,0}(\Om))^n$ and
define the closed sesquilinear form
\[
 \mathfrak{a}:\VV \times\VV\to\CC, \qquad
 \mathfrak{a}(u,v):=\int_\Om \overline{\nabla v}\cdot\nabla u\,dx.
\]
The operator $A$ associated with $\mathfrak{a}$ is the 
\emph{Stokes operator on $\Om$} (with Dirichlet boundary conditions).
It is well-known that $D(A^\einhalb)=\VV$ with equivalent norms
(see \cite{Kato:perturbation}; for the definition of fractional 
domain spaces see Section~\ref{sec:abstract-kato}).
Hence $\dot{\VV}:=\dot{\VV}_2=(\VV,\norm{ \nabla\cdot }_2)^\sim$ equals the 
homogeneous space $\dot{D}(A^\einhalb)$ and the dual space 
$\dot{\VV}^{-1}_2:=(\dot{\VV}_2)'$ can be identified with the 
homogeneous space $(L^2_\si(\Om),\norm{ A^{-\einhalb}\cdot }_2)^\sim$.
Observe that, by Lax--Milgram, a suitable extension $\widetilde{A}$ of 
the operator $A$ acts as an isomorphism $\dot{\VV_2}\to\dot{\VV}^{-1}_2$.
The operator $-A$ generates the bounded analytic semigroup 
$(T(t))=(e^{-tA})$ in $L^2_\si(\Om)$, the \emph{Stokes semigroup}.

\subsection*{$L^q$-theory.}
If there is $q_0\in(2,\infty)$ such that the Helmholtz projection $\PP_{q_0}$
is bounded in $L^{q_0}(\Om)^n$ and there is a bounded analytic semigroup
$T_{q_0}(\cdot)$ in $L^{q_0}$ which is consistent with the Stokes semigroup
in the sense that
\[
 T_{q_0}(t)f=T(t)f, \qquad \mbox{for all $f\in L^2_\si(\Om)\cap L^{q_0}(\Om)$},
\]
then $T_{q_0}(\cdot)$ is called \emph{Stokes semigroup in $L^{q_0}_\si(\Om)$}
or simply \emph{in $L^{q_0}$} and its negative generator $A_{q_0}$ is called 
the \emph{Stokes operator in $L^{q_0}$}. Observe that by interpolation and 
self-duality of the Stokes semigroup we then obtain for any $q\in[q_0',q_0]$
that the Helmholtz projection is $L^q$--bounded and that the Stokes semigroup
extends to a bounded analytic semigroup in $L^q_\si(\Om)$.

\section{Abstract Kato method}\label{sec:abstract-kato}
 
\subsection*{Sectorial operators}
For $0 < \omega \le \pi$ we denote by
\[  \Sec{\om} := \{ z = r e^{i\phi} \suchthat r > 0,\,\, |\phi| <
\omega\}
\]
the open sector of angle $2\omega$ in the complex plane, symmetric about
the positive real axis. In addition we define $\Sec{0} := (0, \infty)$. 
Let $A$ be linear operator on a Banach space $X$. The resolvent set of 
$A$ is denoted by $\varrho(A)$ and its spectrum by $\si(A)$. The operator 
$A$ is called \emph{sectorial of type $\om$}, if 
$\si(A) \subseteq \overline{\Sec{\om}}$ and 
if for all $\nu\in(\om,\pi)$ there is a constant $M$ with 
$\norm{\la (\la+A)^{-1}}\le M$ 
for all $\la\in\Sec{\pi{-}\nu}$. The infimum of all
such angles $\om$ is referred to as the \emph{sectoriality angle of $A$}.

\medskip 

\subsection*{Inter- and Extrapolation spaces}
Given a sectorial operator $A$ on a Banach space $X$, 
for each $n\in\NN$, the space $X_n:=(\DOMAIN(A^n),\norm{(I+A)^n\cdot}_X)$ is a 
Banach space. There are other scales of inter-- and extrapolation spaces.
We give the definitions we need in the sequel, resorting to a
construction in \cite{HaakHaaseKunstmann}: 
Let $A$ be an injective sectorial operator on $X$.
As above, endow $\DOMAIN(A^k)$ with the norm $\norm{(I+A)^{k} \cdot}$ and
$\RANGE(A^k)$ with the corresponding norm $\norm{(I+A^{-1})^{k} \cdot}$.
Let $L := A(I+A)^{-2}$. Then $L(X)=\DOMAIN(A)\cap\RANGE(A)$ and the sum norm
on $\DOMAIN(A)\cap\RANGE(A)$ is equivalent to the norm $\norm{L^{-1}\cdot}=
\norm{(2+A+A^{-1})\cdot}$. Endowing $\XX_1:=\DOMAIN(A)\cap\RANGE(A)$ with this 
norm and letting $\XX_0:=X$ makes $L:\XX_0\to\XX_1$ an isometric 
isomorphism.
Hence, by abstract nonsense we can construct a Banach space $\XX_{-1}$
and an embedding $\iota: \XX_0 \to \XX_{-1}$ together with an 
isometric isomorphism $L_{-1}: \XX_{-1} \to \XX_0$
making the diagram
\[ \xymatrix{
\XX_{-1} \ar@{-->}[r]^{L_{-1}}& **[r]\XX_0=X\\
**[l]X=\XX_0\ar@{->}[u]_{\iota}\ar@{->}[r]^L & **[r]\XX_1=\DOMAIN(A)\cap\RANGE(A)\ar@{-}[u]
}
\]
commute. 
Identifying $\XX_0$ and $\iota\XX_0$ we regard $L$ as a restriction of 
$L_{-1}$. The operator $\AA_{-1}:=L_{-1}^{-1}AL_{-1}$ is an extension of $A$
and again injective and sectorial of the same type in $\XX_{-1}$.

We define recursively for $k\in\NN$ spaces $\XX_{-k}$ and injective 
sectorial operators $\AA_{-k}$ in $\XX_{-k}$, and obtain isometric 
isomorphisms $L_{-k}: \XX_{-k} \to \XX_{-k+1}$, $k\geq 1$.  
In this framework, we now define homogeneous inter-- and 
extrapolation spaces for $k\in\NN$: let 
$\dot{X}_{k} := \AA_{-k}^{-k}(X)$ and 
$\dot{X}_{-k}:= \AA_{-k}^{k}(X)$ with the natural induced norms.
Then we have a scale of spaces
\[
\ldots \emb \dot{X}_n \emb \ldots \emb \dot{X}_1 \emb \dot{X}_0:=X 
\emb \dot{X}_{-1} \emb \ldots \emb \dot{X}_{-n} \emb \ldots
\]
where, for each $n\in\ZZ$, a suitable restriction of 
$\AA_{n}$
acts as an isometric isomorphism $\dot{X}_{n+1} \to \dot{X}_n$.
For each $k\in\NN$, we also let $X_{-k}:=(I+\AA_{-k})^k(X)$ with natural 
norm and denote by $A_{-k}$ the part of $\AA_{-k}$ in $X_{-k}$.
This gives rise to a scale
\[
\ldots \emb {X}_n \emb \ldots \emb {X}_1 \emb {X}_0:=X 
\emb {X}_{-1} \emb \ldots \emb {X}_{-n} \emb \ldots
\]
where, for each $n\in\ZZ$, the operator 
$I+A_{n}$
acts as an isometric isomorphism ${X}_{n+1} \to {X}_n$.

Consequently, if 
the injective operator 
$-A$ generates a semigroup $T(\cdot)$  on $X$ it 
extends in a natural way to a semigroup $T_{-1}(t)=(I{+}A_{-1})\, T(t) \,
(I{+}A_{-1})^{-1}$ on $X_{-1}$, see 
e.g.\cite[Chapter II.5]{EngelNagel} for details.

Notice that $X_{-k} = X + \dot{X}_{-k}$ and $X_{k} = X \cap \dot{X}_{k}$, 
$k\in\NN$.
Moreover, $\dot{X}_k + \dot{X}_{-k} = \XX_{-k}$ and
$\dot{X}_k \cap \dot{X}_{-k} = \DOMAIN(A^k) \cap \RANGE(A^k) =: \XX_{k}$,
$k\in\NN$ (see \cite{HaakHaaseKunstmann,Haase:Buch} for more details).  
We remark that $\dot{X}_k = X_k$ for all $k\in\ZZ$ if $0\in\varrho(A)$.
In any case we have $\norm{x}_{\dot{X}_k} = \norm{A^k x}_X$ 
for $x\in \DOMAIN(A^k)$ and $\norm{x}_{\dot{X}_{-k}} = \norm{A^{-k} x}_X$
for $x\in\RANGE(A^k)$.

Finally we mention that, if $A$ is densely defined with dense range,
we can define the spaces above by completion, i.e. 
\[
\XX_{-k}:=(X,\norm{A^k(I+A)^{-2k}\cdot})^\sim, \qquad 
\dot{X}_k := (\DOMAIN(A^k), \norm{A^k \cdot})^\sim,
\]
\[
\dot{X}_{-k} := (\RANGE(A^k), \norm{A^{-k} \cdot})^\sim,\qquad
X_{-k}:=(X,\norm{(I+A)^{-k}\cdot})^\sim,
\]
for each $k\in\NN$, (see \cite{KaltonKunstmannWeis,
           KaltonWeis:euclidian-structures,
           KunstmannWeis:Levico}).
In this case, we shall also use the notation $\dot{D}(A)$ in place of $\dot{X}_1$ to make
clear with respect to which operator the homogeneous domain space is
taken.

\bigskip

\subsection*{Abstract Kato method}
Let $X$, $Z$, $W$ be Banach spaces and let $\tau \in (0,\infty]$.
Let $-A$ generate a (not necessarily strongly continuous) bounded 
analytic semigroup $T(\cdot)$ on $X$.

Let $B \in B(W, X_{-1})$ and 
$C: X\to Z$ be a closed linear operator that is bounded  $X_1 \to Z$.
Finally let $F: Z \times Z \to W$ be a bilinear map satisfying
$\norm{ F(y, \widetilde y)} \leq K \norm{y}\,\norm{\widetilde y}$ for some
$K>0$. We consider the abstract problem
\begin{equation}\label{eq:kato-abstrakt}
  \left.
  \begin{array}{lcl}
    x'(t) + A x(t) &=& B u(t), \qquad \qquad \; t>0, \\
    x(0) &=& x_0, \\
    y(t) &=& C x(t), \qquad \qquad \; t>0 \\
    u(t) &=& F( y(t), y(t) ) \qquad t>0
  \end{array}
  \right\}
\end{equation}
We seek for mild solutions $x(\cdot)$ in the space $C([0,\tau), X)$,
i.e. for functions $x$ satisfying 
\begin{equation}\label{eq:kato-abstrakt-FP-gleichung}
x(t) = T(t) x_0 + \int_0^t T(t{-}s)B F( Cx(s), Cx(s) )\,ds.
\end{equation}
We shall use the notation $L^p_\al((0,\tau), X) := \{ f \text{
  measurable} \suchthat t^\al f(t) \in L^p( (0,\tau), X)\}$. When $X=\CC$ we also
write $L^p_\al(0,\tau)$.
\begin{theorem}\label{thm:kato-abstrakt}
Let $\tau \in (0,\infty]$ and $p\in (2, \infty]$. Let $\al \ge 0$ such that
$\al+\sfrac1p \in(0, \einhalb)$. We assume
     \let\ALTLABELENUMI\labelenumi \let\ALTTHEENUMI\theenumi
     \renewcommand{\labelenumi}{[A\arabic{enumi}{]}}
     \renewcommand{\theenumi}{[A\arabic{enumi}{]}}
\begin{enumerate}
\item \label{item:C-zul}
     The map $x \mapsto CT(\cdot)x$ is bounded $X \to
     L^p_\al((0,\tau), Z)$.
\item \label{item:B-zul}
     The map $(T_{-1}(\cdot)B) \ast$ is bounded
     $L^{\sfrac{p}2}_{2\al}((0,\tau), W) \to L^\infty((0,\tau), X)$.
\item \label{item:faltung-bdd}
     The map $(C T_{-1}(\cdot)B) \ast$ is bounded $L^{\sfrac{p}2}_{2\al}((0,\tau), W) \to
     L^p_\al((0,\tau), Z)$.
\end{enumerate}
     \let\labelenumi\ALTLABELENUMI
     \let\theenumi\ALTTHEENUMI
Then, under the above assumptions on the operators $B$, $C$ and $F$,
for any initial value $x_0 \in X^\flat := \overline{\DOMAIN(A)}$ 
(the closure being taken in $X$) there exits $\eta \in (0,\tau]$ such that
the abstract problem \eqref{eq:kato-abstrakt} has a unique local mild
solution $x$ in $C([0,\eta), X^\flat)$ satisfying  $Cx \in L^p_\al((0,\eta),Z)$. 
Moreover, if $\norm{x_0}_X$ is sufficiently small, then the solution exists globally.
\end{theorem}

An essential ingredient for the proof is the following lemma 
taking care of the non-linearity (see e.g. \cite[Lemma 1.2.6]{Cannone:Buch}).
\begin{lemma}\label{lem:1-2-6}
  Let $E$ be a Banach space and $\BB: E \times E \to E$ a bilinear map
  with $\norm{ B(e_1,e_2) } \leq \eta \norm{e_1}\, \norm{e_2}$ for all
  $e_1, e_2\in E$. Then, for all $y \in E$ with $\norm{y} <
  \tfrac{1}{4\eta}$ there exists $z\in E$ verifying
  $z = y + \BB(z,z)$ and $\norm{z}_E \leq 2 \norm{y}_E$.
\end{lemma}
The lemma is shown by resorting to Banach's fixed-point theorem on a small ball within $E$.
\begin{proof}[Proof of {Theorem~\ref{thm:kato-abstrakt}}]
{\em (Existence)}
Let $\eta>0$ and $E := L^p_\al( (0,\eta), Z)$ and consider the
bilinear map
\[
\BB := \left\{\begin{array}{lcl}
    E \times E & \to & E \cr
    (x, \widetilde x) &\mapsto & (C T_{-1}(\cdot)B) \ast F(x,\widetilde x).
 \end{array}\right.
\]
For $x_0 \in X$ we have $y := CT(\cdot)x_0 \in E$ by \ref{item:C-zul}. 
For $x,\widetilde x\in E$ we have $F(x,\widetilde x)\in L^{\sfrac{p}2}_{2\al}((0,\eta),W)$.
Moreover, $\BB$ is bounded by \ref{item:faltung-bdd}. If $p<\infty$,

\[
\norm{y}_E = \biggl(\int_0^{\eta} \norm{ t^\al CT(t)x_0 }^p_Z \,dt\biggr)^\sfrac{1}p
\]
becomes small if $\eta>0$ is small enough. 
In the case that $p=\infty$, notice that 
by \ref{item:C-zul}, $t^\al T(t)x_0$ is bounded in $Z$ and that for $x_0
\in \DOMAIN(A)$ we have $t^\al T(t) x_0 \to 0$ in $X_1$ for $t\to 0+$ 
since $AT(t) x_0$ is bounded near the origin and $\al >0$.
Thus, for all  $p \in [2,\infty]$ and $x_0 \in X^\flat$, we can 
make $\norm{y}_E$ 
arbitrarily 
small choosing $\eta>0$ small enough. 
If, on the other hand, $\tau=\infty$ and $\norm{x_0}$ is small enough, 
assumption \ref{item:C-zul} allows to take $\eta=\infty$.

\medskip

In any case Lemma~\ref{lem:1-2-6} applies and shows existence of a solution
$z\in E$ satisfying $z = y + \BB(z,z)$. Now put
\[
x(t) := T(t) x_0 + \int_0^t T(t{-}s) B F(z(s), z(s) )\,ds
\]
Then, by \ref{item:B-zul} $x \in L^\infty((0,\eta), X)$.
By definition of $x$ and the fixed--point equation satisfied by $z$,
$z(t) = C x(t)$  (recall that $C$ is closed as operator
$X \to Z$). Thus, $x(\cdot)$ is a mild solution of the abstract problem
\eqref{thm:kato-abstrakt} as claimed.

\smallskip

{\em (Continuity)}
To see that $x$ is continuous with values in $X^\flat$ we go again
through the fixed-point argument. We employ the following ad-hoc
notation: for a Banach space $Y$ let
\[
C_0([0,\tau), Y) = \{ y: [0,\tau)\to Y \suchthat y \in C_b([0,\tau), Y) \text{ and } y(0)=0 \},
\]
endowed with the supremum norm, and let  
\[
C_{0,\al}([0,\tau), Y) = \{ y: [0,\tau)\to Y \suchthat t^\al y(t) \in  C_0([0,\tau), Y) \},
\]
endowed with the weighted supremum norm $\norm{y}_{C_{0,\al}} =
\sup_{t\in [0,\tau)} \norm{t^\al y(t) }_Y$.  First notice that for
$x_0 \in X^\flat$, the map $t \mapsto CT(t)x_0$ defines an element of
$C_{0,\al}([0,\tau), Z)$. The norm estimate is clear by
\ref{item:C-zul}.  Strong continuity of the semigroup implies that for
$x_0 \in \DOMAIN(A^\flat)$, the trajectory $T(t)x_0$ is bounded and
continuous within $[\DOMAIN(A^\flat)]$ and so $t^\al T(t) x_0$ defines
an element of $C_0([0,\min(\tau, r)), [\DOMAIN(A^\flat)])$ for all
$r>0$.  Since $C\in B(X_1, Z)$ the first claim follows by letting 
$r\to\infty$ and using the density of $\DOMAIN(A^\flat)$ in $X^\flat$.
Next, observe that the bilinear map $(u, v) \mapsto F(u,v)$ defines a
continuous map from $C_{0,\al}([0,\tau), Z) \times C_{0,\al}([0,\tau),Z)$ 
to $C_{0,2\al}([0,\tau), W)$.  Finally, we show that
\[
\begin{split}
  T(\cdot)B \ast: & \quad C_{0,2\al}([0,\tau), W) \to C_0([0,\tau), X) \\
 CT(\cdot)B \ast: & \quad C_{0,2\al}([0,\tau), W) \to C_{0,\al}([0,\tau), Z).
\end{split}
\]
Since continuity is a local property we may assume $\tau<\infty$ without 
loss of generality. 
By \cite[Theorem 1.8, 1.9]{HaakKunstmann:weighted-Lp-admiss},
\ref{item:B-zul} implies $\RANGE(B) \subseteq (\dot X_{-1},
X)_{\theta,\infty}=:\widetilde W$ with $\theta=\sfrac2p+2\al \in
(0,1)$. Let $u \in  C_{0,2\al}([0,\tau), W)$. Then
$Bu \in C_{0,2\al}([0,\tau), \widetilde W)$ and so 
analyticity of the semigroup implies
$\norm{ T(t-s) B u(s) }_X \le C (t{-}s)^{\theta-1} s^{-2\al}$, 
the integral is thus absolutely convergent within $X$. Again
by analyticity, $T(t{-}s)B u(s) \in \DOMAIN(A)$ and since 
$X^\flat$ is a closed subspace of $X$, 
$(T(\cdot)B\ast)(t) u \in X^\flat$ for every $t\in [0,\tau)$.
Now,
\begin{eqnarray*}
&&   \Bignorm{\bigl(T(\cdot)B\ast u\bigr)(t) - \bigl(T(\cdot)B\ast u\bigr)(r)}_{X^\flat} \\
&\le&  \Bignorm{ \bigl(T(t{-}r)-I\bigr) (T(\cdot)B\ast u)(r) }_{X^\flat}
  +  M \int_r^t (t{-}s)^{\theta-1} s^{-2\al}\,ds
\end{eqnarray*}
for all $t>r$, and so strong continuity of the semigroup on $X^\flat$
shows left continuity. A similar argument yields right continuity and 
so the first claim follows. The second is shown similarly: by \cite[Theorem
1.7, 1.9]{HaakKunstmann:weighted-Lp-admiss}, \ref{item:C-zul} implies
that $C$ is bounded $\widetilde Z \to Z$ where $\widetilde Z = (X,
\dot X_1)_{\al+\sfrac1p, 1}$. Notice that 
$\norm{T(t)}_{\widetilde W  \to \widetilde Z} \le c t^{\sfrac1p+\al-1}$ for 
$t>0$ which implies that the convolution is an absolutely convergent 
integral within $Z$. Moreover, for small $\eps>0$ ,
$T(\cdot)B\ast u$ is absolutely convergent within $\widetilde W_{1-\eps}$
and so
\begin{eqnarray*}
 &&  \Bignorm{t^\al \bigl(T(\cdot)B\ast u\bigr)(t) - r^\al \bigl(T(\cdot)B\ast u\bigr)(r)}_{\widetilde Z} \\
&\le& \Bignorm{ \bigl(t^\al T(t{-}r)-r^\al\bigr) (T(\cdot)B\ast u)(r) }_{\widetilde Z}
  +  M t^\al \int_r^t (t{-}s)^{\al+\sfrac1p-1} s^{-2\al}\,ds 
\end{eqnarray*}
Since $\norm{AT(t)}_{\widetilde W_{1-\eps} \to \widetilde Z} \le C
t^{\al+\sfrac1p-1-\eps}$ is integrable at the origin, left continuity
follows.  A similar argument for right continuity shows the second
claim. We conclude by the fixed point equation
(\ref{eq:kato-abstrakt-FP-gleichung}) that the solution $x$ is
continuous in $X^\flat$.

\smallskip

{\em (Uniqueness)}
Assume the existence of two solutions $u,v$ to (\ref{eq:kato-abstrakt}) in 
$C([0,\eta), X^\flat)$ satisfying both $Cu, Cv \in L^p_\al((0,\eta), Z)$ and therefore
satisfying both the fixed point equation (\ref{eq:kato-abstrakt-FP-gleichung}). 
Let $\eta_0\in (0,\eta)$. Using bilinearity and continuity of $F$ and assumption 
\ref{item:faltung-bdd} we obtain
\begin{eqnarray*}
&&      \norm{C(u-v)}_{L^{p}_{\al}((0,\eta_0), Z)} \\
& = & \norm{ (C T(\cdot)B) \ast (F(Cu,Cu)- F(Cv,Cv)) }_{L^{p}_{\al}((0,\eta_0), Z)} \\
&\le& \norm{ (C T(\cdot)B) \ast (F(Cu,C(u-v))) }_{L^{p}_{\al}((0,\eta_0), Z)} \\
       && + \; \norm{ (C T(\cdot)B) \ast (F(C(u-v),Cu)) }_{L^{p}_{\al}((0,\eta_0), Z)} \\
&\le& M \bigl( \norm{  F(Cu,C(u-v)) }_{L^{\sfrac{p}2}_{2\al}((0,\eta_0), W)} 
       + \; \norm{ F(C(u-v),Cu) }_{L^{\sfrac{p}2}_{2\al}((0,\eta_0), W)} \bigr) \\
&\le& M \norm{F} \; \bigl( \norm{Cu}_{L^{p}_{\al}((0,\eta_0), Z)} \norm{C(u-v)}_{L^{p}_{\al}((0,\eta_0), Z)} \\
       && \qquad \;\; + \; \norm{Cv}_{L^{p}_{\al}((0,\eta_0), Z)} \norm{C(u-v)}_{L^{p}_{\al}((0,\eta_0), Z)} \bigr)\\
& = & M \norm{F} \; \bigl( \norm{Cu}_{L^{p}_{\al}((0,\eta_0), Z)} + \norm{Cv}_{L^{p}_{\al}((0,\eta_0), Z)} \bigr) \; 
         \norm{C(u-v)}_{L^{p}_{\al}((0,\eta_0), Z)}.
\end{eqnarray*}
Now choosing $\eta_0>0$ small enough makes $\norm{Cu}_{L^{p}_{\al}((0,\eta_0), Z)}$ and 
$\norm{Cv}_{L^{p}_{\al}((0,\eta_0), Z)}$ 
arbitrarily 
small which allows to conclude 
$Cu=Cv$ in $L^p_\al((0,\eta_0), Z)$. 
For $p<\infty$ this smallness is immediate; for $p=\infty$ we argue with
$u-T(\cdot)x_0\in C_0([0,\eta),X^\flat)$
 as above. Thus,
\begin{eqnarray*}
u(t) &=&  T(t) x_0 + \int_0^t T(t{-}s)B F( Cu(s), Cu(s) )\,ds \\
     &=&  T(t) x_0 + \int_0^t T(t{-}s)B F( Cv(s), Cv(s) )\,ds
      =   v(t)
\end{eqnarray*}
for $t\in [0,\eta_0)$. Repeating the argument with $x_0 := u(\eta_0) \in X^\flat$ yields 
uniqueness of the solution as claimed.
\end{proof}

\begin{remark}
Notice that in a setting of linear systems theory, assumptions
\ref{item:C-zul} and \ref{item:B-zul} of the theorem mean 
weighted admissibility conditions for the observation operator $C$ and
the control operator $B$. We refer to 
\cite{HaakKunstmann:weighted-Lp-admiss} for more details. 
\end{remark}
\begin{remark}\label{rem:weissler}
In the applications to \eqref{eq:nse}, the operators $C$ and $B$ are 
suitable identity operators. Weissler's result \cite{Weissler:NSE}
assumes continuity of the bilinearity $Z\times X\to W$.  

Observe that, for general operators $B$ and $C$, this leads to a 
different setting, whereas we are working entirely with the space 
$Z$ for the fixed point argument. In applications to \eqref{eq:nse} 
this has the advantage that (tensor) products need only be defined
for elements of $Z$, and that we can allow for spaces $X$ with very
rough initial data (see Section~\ref{sec:applications}).
\end{remark}
\begin{corollary}\label{cor:Kato-Abstrakt-mit-mehr-Kontrolltermen}
Let additionally to the situation in Theorem~\ref{thm:kato-abstrakt} 
Banach spaces $W^{(1)}, \ldots, W^{(m)}$ and operators 
$B_j \in B(W^{(j)}, X_{-1})$ be given and consider the abstract problem
\begin{equation}\label{eq:Kato-inhomogen}
  \left.
  \begin{array}{lcl}
    x'(t) + A x(t) &=& B u(t) + \sum\limits_{j=1}^m B_j f_j(t), \qquad t>0, \\[-1.5ex]
    x(0) &=& x_0, \\
    y(t) &=& C x(t), \qquad \qquad \qquad \qquad \; t>0\\
    u(t) &=& F( y(t), y(t) ) \qquad \qquad \qquad t>0
  \end{array}
  \right\}
\end{equation}
where the 'inhomogeneities' $f_j$ satisfy 
$f_j\in L^{p_j}_{\beta_j}((0,\tau), W^{(j)})$
for some $p_j, \beta_j$ with $\sfrac{1}{p_j} + \beta_j \in (0,1)$. 
Moreover we require 
     \let\ALTLABELENUMI\labelenumi \let\ALTTHEENUMI\theenumi
     \renewcommand{\labelenumi}{[$A_j$\arabic{enumi}{]}}
     \renewcommand{\theenumi}{[$A_j$\arabic{enumi}{]}}
\begin{enumerate} \setcounter{enumi}{1}
 \item \label{item:A2-inhomogen}
      The maps $(T_{-1}(\cdot)B_j) \ast$ are bounded
      $L^{p_j}_{\beta_j}((0,\tau), W^{(j)}) \to L^\infty((0,\tau), X)$.
\item \label{item:A3-inhomogen}
     The maps $(CT(\cdot)B_j) \ast: L^{p_j}_{\beta_j}((0,\tau), W^{(j)}) 
          \to L^p_\al((0,\tau),Z)$  are bounded.
\end{enumerate}
     \let\labelenumi\ALTLABELENUMI
     \let\theenumi\ALTTHEENUMI
for all $j=1,\ldots, m$.
Then time-local mild solutions always exist in case $p<\infty$. In case 
$p=\infty$ or in order to obtain global solutions a smallness condition 
on the norms of the functions $f_j$ has to be imposed ($j=1,\ldots, m$).
\end{corollary}
\begin{proof}
As for Theorem~\ref{thm:kato-abstrakt} but with 
$y = CT(\cdot)x_0 + \sum\limits_{j=1}^m (CT(\cdot)B_j) \ast f_j$.
\end{proof}

Before coming to applications in Section~\ref{sec:applications} we sum up
necessary and sufficient conditions for  the assumptions 
\ref{item:C-zul} -- \ref{item:faltung-bdd} in Theorem~\ref{thm:kato-abstrakt}.
Notice that \ref{item:A2-inhomogen} and \ref{item:B-zul} are of the same 
type whereas \ref{item:faltung-bdd} is a special case of
\ref{item:A3-inhomogen}. 
Throughout the rest of this article, 
for any interpolation couple $(E, F)$
of Banach spaces we denote  by  $(E,F)_{\si, p}$ the real
interpolation space between $E$ and $F$.
For references and
details on the real interpolation method see e.g. \cite{BerghLoefstroem,Lunardi:Buch,Triebel:interpolation}.

\begin{theorem}\label{thm:discussion-A1-A3}
Let $p \in(2,\infty]$ and $\al \ge 0$ such that $\al+\sfrac{1}p \in (0,\einhalb)$. 
 \begin{enumerate}
  \item \label{item:discussion-A1}
    If {\rm \ref{item:C-zul}} holds for $\tau=\infty$, then $C$ is bounded in
    the norm $(X, \dot{X}_1)_{\al+\sfrac{1}p, 1} \to Z$. 
    The converse is true provided that $X \emb (\dot{X}_{-1},
    \dot{X}_1)_{\einhalb, p}$.
  \item  \label{item:discussion-A2}
    If {\rm \ref{item:B-zul}} holds for $\tau=\infty$, then $B$ is bounded in
    the norm $W \to (\dot{X}_{-1}, X)_{2(\al+\sfrac{1}p), \infty}$.
    The converse is true in case $\al>0$ or in case $\al=0$ and $(\dot{X}_{-1},
    \dot{X}_1)_{\einhalb, \sfrac{p}2} \emb X$.
  \item  \label{item:discussion-A3j}
    The map $(CT(\cdot)B) \ast: L^q_\beta((0,\tau), W) \to L^p_\al((0,\tau), Z)$
    is bounded 
    provided that $\norm{C T(t)B}_{W\to Z} \leq c \,t^{-\ga}$ for $t\in(0,\tau)$
    and that $\beta+\ga+\sfrac1q = 1 +\al+\sfrac{1}p$ where 
    $\ga\in (0,1)$ and $\al,\beta >0$.
 \end{enumerate}
\end{theorem}

A proof of the theorem and some additional results will be provided in 
Section~\ref{sec:proof-of-kato-discussion}. As the proof will actually show, 
the restriction $p>2$ (instead of $p>1$) and $\al+\sfrac1p<\einhalb$ 
(instead of $<1$) in the above formulation is only due to 
the bilinear structure which forces to consider the parameters 
$\sfrac{p}2$ and $2\al$ in part \ref{item:discussion-A2}. We mention 
that in part \ref{item:discussion-A1} (and in part \ref{item:discussion-A2} 
in case $\al=0$) of the 
theorem 
the embedding assumption on the space $X$ is optimal. This follows by
choosing $C=A^{\al+\sfrac1p}$ and $B=A^{1-\sfrac2p}$, see also the
discussion in \cite[Section 1]{HaakKunstmann:weighted-Lp-admiss}.

\begin{remark}\label{rem:discussion-A3}
As mentioned above, condition \ref{item:A3-inhomogen} contains 
condition \ref{item:faltung-bdd} as special cases letting 
$p_j=\sfrac{p}2$ and $\beta_j=2\al$. Here, and for one direction of 
part \ref{item:discussion-A2} in case $\al>0$, the proof 
is based 
on a classical one-dimensional convolution estimate due to {\sc Hardy} 
and {\sc Littlewood}, see Lemma~\ref{lem:HardyLittlewood}.
\end{remark}

\begin{remark}\label{rem:discussion-finite-time}
  Supposing that $\tau<\infty$ in
  Theorem~\ref{thm:discussion-A1-A3}~\ref{item:discussion-A3j} and
  $\beta+\ga+\sfrac1q \le 1 +\al+\sfrac{1}p$, one finds $\widetilde
  \ga>\ga$ such that $\norm{CT(t)B}_{W\to Z} \le \widetilde
  c\,t^{-\widetilde \ga}$ and $\beta+\widetilde \ga+\sfrac1q = 1
  +\al+\sfrac{1}p$.  Concerning the parts \ref{item:discussion-A1}
  and~\ref{item:discussion-A2} of Theorem~\ref{thm:discussion-A1-A3} we
  remark that if \ref{item:C-zul} or \ref{item:B-zul} hold on
  $(0,\tau)$ for $A$, then they hold for any $\nu>0$ for $\nu{+}A$.
  From this it is clear how to modify the resolvent conditions in
  section~\ref{sec:proof-of-kato-discussion}: the homogeneous spaces
  $\dot X_{-1}$ and $\dot X_{1}$ have to be replaced by their
  inhomogeneous counterparts $X_{-1}$ and $X_{1}$ (see also
  \cite{HaakHaaseKunstmann}). In particular, if $\tau<\infty$, one can
  without loss of generality assume $A$ to be boundedly invertible
  (see also \cite[Lemma 1.3]{HaakKunstmann:weighted-Lp-admiss}).
\end{remark}

\subsection*{The r\^ole of maximal regularity for recovering the pressure terms}

From now on we shall always use $C=\text{Id}_Z$ and $B=\text{Id}_W$, i.e.
we suppose $X_1 \emb Z$ and $W \emb X_{-1}$. Consequently, the semigroup 
$T(t)$ acts (pointwise as a bounded operator) $X\to Z$ and $W\to X$. In 
this setting, the abstract Cauchy problem \eqref{eq:kato-abstrakt} takes 
the form
\begin{equation} \label{eq:abstract-bilinear}
\left.
  \begin{array}{lcl}
      x'(t) + A x(t) & = &F(x(t), x(t)) \cr
      x(0) & = &x_0      
  \end{array}
\right\}
\end{equation}
\begin{proposition}\label{prop:simplification-A1-A2}
The embedding assumptions in Theorem~\ref{thm:discussion-A1-A3} \ref{item:discussion-A1} 
and \ref{item:discussion-A2} concerning the spaces $Z$ and $W$ are equivalent to 
pointwise growth estimates for the semigroup. Indeed, for $\si,\theta\in (0,1)$ one has
\[
(X, \dot{X}_1)_{\si, 1} \emb Z
\quad \text{if and only if}\quad
\norm{T(t)}_{X\to Z} \le c \,t^{-\si}
\]
and
\[
W \emb (\dot{X}_{-1}, X)_{\theta, \infty} 
\quad \text{if and only if}\quad
\norm{T(t)}_{W \to X} \le c \,t^{\theta-1}
\]
\end{proposition} 
\begin{proof}
Since the semigroup is bounded and analytic, one has the elementary estimates
\[
\norm{T(t)}_{\dot X_n \to \dot X_{n+1}} \le C \,t^{-1}
\quad\text{and}\quad
\norm{T(t)}_{\dot X_n \to \dot X_n} \le M
\]
for $n\in\ZZ$. Therefore, 
the embedding properties for $Z$ and $W$ imply the growth estimates of the semigroup 
acting $X\to Z$ and $W\to X$ by interpolation.
Conversely, the estimate $\norm{T(t)}_{X\to Z} \le c \,t^{-\si}$ for $t>0$ implies
\begin{align*}
\norm{x}_Z 
&  =  \widetilde c \; \biggnorm{\int_0^\infty A T(2t)x \,dt}_Z
  \le \widetilde c \int_0^\infty \bignorm{ A T(2t)x}_Z \,dt \\
& \le c \int_0^\infty t^{-\si} \bignorm{ A T(t)x}_X \,dt
= c \;\norm{x}_{(X,\dot X_1)_{\si,1}}
\end{align*}
for $x \in X_1 = \dot X_1 \cap X$ which is dense in $(X,\dot X_1)_{\si,1}$.
Finally, the estimate $\norm{T(t)}_{W\to X} \le c \,t^{\theta-1}$ for $t>0$ implies
\[
\norm{x}_{(\dot X_{-1}, X)_{\theta,\infty}} 
= \sup_{t>0} \; t^{-\theta} \norm{t A T(t) x}_{\dot X_{-1}} 
= \sup_{t>0} \; t^{1-\theta} \norm{T(t) x}_{X} \le c \, \norm{x}_W
\]
which finishes the proof.
\end{proof}

Given an abstract Cauchy problem of the form
\begin{equation}\label{eq:ACP}
    x'(t) + A x(t) = f(t),\qquad x(0)=0
\end{equation}
on a Banach space $W$, we say that $A$ has {\em maximal $L^p$--regularity},
$p\in [1,\infty]$ if the mild solution $x$ to \eqref{eq:ACP} satisfies 
$x', Ax \in L^p((0,\tau),W)$ whenever $f \in L^p((0,\tau),W)$. 
We refer to \cite{Amann:Quasilin-Buch,DaPratoGrisvard,Dore:LectureNotes,
Haase:Buch,KunstmannWeis:Levico,Weis:FM} for this relation, the problem of 
maximal regularity, characterisation results, and further references on the 
subject.

\begin{theorem}\label{thm:strong-solutions}
Suppose $\tau\in (0,\infty]$, $p \in(2,\infty]$ and 
$\al \ge 0$ 
such that  $\al+\sfrac{1}p \in (0,\einhalb)$.  Let $X, Z, W$ be Banach spaces satisfying
\begin{eqnarray}
 W &\emb& (\dot{X}_{-1}, X)_{2(\al+\sfrac1p), \infty},             \label{eq:assumption-W} \\
(X, \dot{X}_1)_{\al+\sfrac{1}p, 1} & \emb&  Z, \qquad \text{and} \label{eq:assumption-Z} \\
 X  &\emb& (\dot{X}_{-1}, \dot{X}_1)_{\einhalb, p}               \label{eq:assumption-X}
\end{eqnarray}
and assume that $-A$ is injective and generates consistent bounded analytic
semigroups on $X$ and $W$. Let $A$ have maximal $L^{\sfrac{p}2}$--regularity on $W$.
Then for every $x \in X^\flat$, the abstract problem \eqref{eq:abstract-bilinear}
has a unique time-local mild solution 
\[
   x \in C([0,\eta), X) \cap L^p_\al((0,\eta), Z) \cap L^{\sfrac{p}2}_{2\al}((0,\eta), W)
\]
that satisfies $x',Ax\in L^{\sfrac{p}2}_{2\al}((0,\eta),W) + L^{p}_{\al+1}((0,\eta),Z)$.

\end{theorem}
\begin{proof}
By Theorem~\ref{thm:discussion-A1-A3}, equation \eqref{eq:abstract-bilinear} has a 
time-local mild solution $x$ as claimed for some $\eta\in(0,\tau)$. 
The {\sc Pr\"u\ss{}-Simonett} theorem (see \cite[Theorem 2.4]{PruessSimonett}, also
\cite[Theorem 1.13]{HaakKunstmann:weighted-Lp-admiss}) shows that maximal 
$L^\sfrac{p}2$--regularity of $A$ in $W$ induces
maximal $L^\sfrac{p}2_{2\al}$--regularity in $W$ 
(recall that $\al+\sfrac1p <\einhalb$); 
this result is also true in case $p=\infty$, as an inspection 
of the proof shows 
(it is actually even easier to prove than for finite $p$). 
Let $x$ be the mild solution to \eqref{eq:abstract-bilinear}. Writing
\[  
    x(t) = T(t)x_0 + \bigl(T\ast F(x,x)\bigr)(t) = x_1(t) + x_2(t) 
\]
one deduces from maximal $L^\sfrac{p}2_{2\al}(W)$--regularity that 
$x_2 \in W_1 = \DOMAIN(A_W)$ a.e. and that $x_2$ is a.e. differentiable 
satisfying $x_2', Ax_2 \in L^\sfrac{p}2_{2\al}((0,\eta),W)$. \\

Using Proposition~\ref{prop:homogene-interpol-T}, we have
$x_1' = -A x_1 \in L^p_{\al+1}((0,\eta), (X, \dot X_1)_{\al+\sfrac1p, 1} )$ 
for 
\[
x_0 \in 
\bigl( (\dot X_{-1}, X)_{\al+\sfrac1p, 1}, (\dot X_1, \dot X_2)_{\al+\sfrac1p, 1} \bigr)_{1-\einhalb(\al+\sfrac1p+1), p}
=  (\dot X_{-1}, \dot X_1)_{\einhalb, p}
\]
where the equality is due to reiteration. 
By \eqref{eq:assumption-X}, this condition 
holds for $x_0\in X^\flat$. Finally, assumption 
\eqref{eq:assumption-Z} finishes the proof.
\end{proof}

Observe if one has  $X \emb (\dot{X}_{-1}, \dot{X}_1)_{\einhalb, \sfrac{p}2}$ 
in place of (\ref{eq:assumption-X}) one obtains 
$x', Ax \in L^{\sfrac{p}2}_{2\al}((0,\eta),W)$ by similar arguments. This applies in 
particular in case $p=\infty$.

\medskip

For $p \in[2,\infty)$, the {\sc da$\,$Prato-Grisvard} theorem 
(\cite{DaPratoGrisvard}, see also \cite[Theorem 9.3.9]{Haase:Buch}) 
provides several function spaces in which negative generators of 
analytic semigroups have maximal $L^{\sfrac{p}2}$--regularity. In our 
situation, a particularly import class are real interpolation spaces 
of the form $(\dot X_{-k}, \dot X_k)_{\theta,r}$ for some $k\in \NN$ and
$\theta \in (0,1)$, $r\in [1,\infty]$. When $1<p<\infty$,
 maximal $L^p$--regularity is independent of $p$ 
(see \cite{BenedekCalderonPanzone}). In case 
$p=\infty$, the following lemma may be used 
to verify $L^\infty$--maximal regularity.
\begin{lemma}\label{lem:L-infty-MR}
Let the injective operator $-A$ generate a (not necessarily densely
defined) bounded analytic semigroup on $W$ and let $U$ be a Banach
space, such that $\norm{T(t)}_{W \to U} \leq c\,t^{-1}$ for some $c>0$.
If $(W, \dot W_2)_{\einhalb, \infty} \emb U$, then
\[
\ess\sup_{t>0} \;\Bignorm{ \int_0^t T(t{-}s) w(s)\,ds}_U \leq C \,\norm{w}_{L^\infty(\RR_+, W)}
\]
for all $w\in L^\infty(\RR_+, W)$.
\end{lemma}
\begin{proof}
It is clearly sufficient to verify
\[
\Bignorm{\int_0^\infty T(s) w(s) \,ds }_{(W, \dot W_2)_{\einhalb,\infty}} 
  \leq C \norm{w}_{L^\infty(\RR_+, W)}
\]
for all $w \in L^\infty(\RR_+, W)$. Using Proposition~\ref{prop:homogene-interpol-T} one has
\begin{eqnarray*}
       \Bignorm{\int_0^\infty T(s) w(s) \,ds }_{(W, \dot W_2)_{\einhalb,\infty}} 
&\sim&  \Bignorm{t \mapsto t A^2 T(t) \int_0^\infty T(s) w(s) \,ds}_{L^\infty(\RR_+, W)} \\
& = &  \Bignorm{t \mapsto \int_0^\infty t A^2 T(t{+}s) w(s) \,ds}_{L^\infty(\RR_+, W)} \\
&\leq&  \ess\sup_{t>0} \int_0^\infty \tfrac{t}{(t+s)^2} \norm{w(s)}_W \,ds  \\
&\leq&  \norm{w}_{L^\infty(\RR_+, W)} \int_0^\infty \tfrac{1}{(1+\si)^2} \,d\si
\end{eqnarray*}
by substituting $s=t \si$.
\end{proof}

\begin{remark}\label{rem:get-pressure-back}

By Theorem~\ref{thm:strong-solutions} one can give a sense to the
differential equation
 in (\ref{eq:ACP})
for a.e. $t>0$ in the time interval under consideration. In applications 
to the Navier-Stokes  
equations (\ref{eq:nse})
this means that
\[
u'(t) + A u(t) + P\nabla\cdot(u(t)\otimes u(t)) - P f(t) = 0
\]
for a.e. $t>0$. Interpreting $A$ as $-P\Delta$ and the operator $P$ as restriction
$P: \sD'(\Om)^n \to \sD_\si'(\Om)$ (see e.g. \cite{Monniaux:arbitrary,Sohr:Buch})
we are led to
\[
P \bigl( u_t -\Delta u(t) + \nabla\cdot (u\otimes u) -f \bigr)=0
\]
if $u_t -\Delta u(t) + \nabla\cdot (u\otimes u) -f \in \sD'(\Om)$ at a 
fixed time $t>0$. Now the pressure term $\nabla p$ can be recovered by 
Theorem~\ref{thm:deRham} which passes from (\ref{eq:nseP}) back to (\ref{eq:nse}).

\medskip

The equality $A=-P\Delta$ is no problem in case $\Om=\RR^n$ since then 
$\Delta$ commutes with $P$. If $\Om\subseteq \RR^n$ is bounded or an exterior domain 
with $\partial \Om \in C^{1,1}$, equality $A=-\PP\Delta_D$ holds on 
$\DOMAIN(A_q) = W_q^2(\Om)^n \cap W_{q,0}^1(\Om)^n \cap L^q_\si(\Om)$ for $q\in (1,\infty)$ 
where $\Delta_D$ denotes the Dirichlet Laplacian on $\Omega$. On arbitrary domains
$\Om \subseteq \RR^3$, equality $A=-P\Delta$  holds on $\VV = \VV_2$, 
see \cite{Monniaux:arbitrary}.
\end{remark}

\section{Application to the Navier-Stokes equations}\label{sec:applications}

In this section we apply the abstract result to the Navier--Stokes 
equations \eqref{eq:nse}, where \eqref{eq:kato-abstrakt} corresponds 
to \eqref{eq:nseP} and \eqref{eq:kato-abstrakt-FP-gleichung} corresponds 
to \eqref{eq:fix-point}. In these applications we always have 
$C=\text{Id}_Z$ and $B=\text{Id}_W$ which means that the necessary
conditions in Theorem~\ref{thm:discussion-A1-A3} 
\ref{item:discussion-A1} and \ref{item:discussion-A2} boil down to
continuous embeddings or via 
Proposition~\ref{prop:simplification-A1-A2} to decay estimates for the semigroup.

It turns out that the choice of the ``auxiliary space'' $Z$
is most significant. The structure of the map $F$ then determines the 
space $W$, and one can calculate the exponent $\ga$ for which
$\norm{T(t)}_{W \to Z} \leq C \,t^{-\ga}$
 holds. 
Depending on the context, 
this may hold on $(0,\infty)$ or on bounded time intervals $(0,\tau)$
where $\tau<\infty$. Observe that an application of  
Theorem~\ref{thm:discussion-A1-A3}~\ref{item:discussion-A3j} 
(and Remark~\ref{rem:discussion-A3}) requires 
$\ga\in (\einhalb, 1)$ and restricts $\al$ and $p$ to 
$\al+\sfrac{1}p\leq 1-\ga$ for local
solutions and to $\al+\sfrac{1}p=1-\ga$ for global solutions. 
Nevertheless, we have some freedom for the choice of $\al$ and $p$. 
Once $\al$ and $p$ are fixed, Theorem~\ref{thm:discussion-A1-A3}
\ref{item:discussion-A1} and \ref{item:discussion-A2} 
allow to adjust the space $X$ for initial values appropriately.

In the sequel we discuss various choices of $Z$ on $\RR^n$ and on
domains.  The common approach covers some known results, provides new
proofs for other known results, but it also yields new results on
$\RR^n$ and on domains. 

\medskip
\subsection{Lebesgue spaces  on $\RR^n$}\label{sec:lebesgue}

Here and in the following subsections we consider the Navier--Stokes
equations \eqref{eq:nse} on $\RR^n$, $n\geq 2$. For simplicity we
shall omit $\RR^n$ and superscripts $n$ or $n\times n$ in notation. 
On $\RR^n$ the Helmholtz projection commutes with the Laplacian
$\Delta$ and is bounded on $L^q$ for $1<q<\infty$.

Let $q\in(n,\infty)$ and consider the case $Z = L^q$. For $u,v\in Z$ we 
then have $\nabla\cdot(u\otimes v) \in \dot{H}^{-1}_{\sfrac{q}2} =: W$. 

Notice that 
$\norm{T(t)}_{\dot{H}_{\sfrac{q}2}^{-1} \to \dot{H}_{\sfrac{q}2}^{\delta}} 
\leq c \,t^{-(1+\delta)/2}$, $t>0$, 
and that $\dot{H}_{\sfrac{q}2}^{\delta} \emb L^q$ provided 
$\sfrac{1}q=\sfrac{2}q-\sfrac{\delta}n$, i.e. provided
$\delta = \sfrac{n}q$ 
(see e.g.  \cite[Theorems 2.7.1 and 5.2.5]{Triebel:FunctionSpaces}).
We obtain $\norm{T(t)}_{W\to Z}\le c \,t^{-\ga}$, $t>0$, where
$\ga=\frac{1}2+\frac{n}{2q}\in(\einhalb,1)$. Hence we should have
$\al+\sfrac{1}p=1-\ga=\frac{1}2-\frac{n}{2q}$ which restricts $p$ to
$p\in[\frac{2q}{q-n},\infty]$. For such a $p$, consider the space 
$X=\dot{B}^{-1+\sfrac{n}q}_{q,p}$. Then
\[
(\dot{X}_{-1}, \dot{X}_{1})_{\einhalb, p} = X
\quad\text{and}\quad
(\dot{X}_{-1}, \dot{X}_{1})_{\einhalb, \sfrac{p}2} =
\dot{B}^{-1+\sfrac{n}q}_{q,\sfrac{p}2} \emb X
\]
 whence Theorem~\ref{thm:discussion-A1-A3} allows to verify \ref{item:C-zul} and
\ref{item:B-zul} of Theorem~\ref{thm:kato-abstrakt} easily. 
Indeed, we have $(X,\dot{X}_1)_{\al+\sfrac{1}p, 1} = \dot{B}^s_{q,1}$ with
$s=2\al+\sfrac{2}p+\sfrac{n}q-1$ and $C=\text{Id}_Z$ certainly satisfies the
condition of Theorem~\ref{thm:discussion-A1-A3} \ref{item:discussion-A1}
if $\dot{B}^s_{q,1}$ embeds into $Z$ which is the case if $s= 0$,
i.e. if  $\al+\sfrac{1}p = \frac12 - \frac{n}{2q}$.

Moreover, 
$(\dot{X}_{-1}, X)_{2(\al+\sfrac{1}p), \infty} = \dot{B}^t_{q,\infty}$
with $t=4(\al+\sfrac{1}p) + \sfrac{n}q - 3$ whence $B$ satisfies the
condition in Theorem~\ref{thm:discussion-A1-A3} \ref{item:discussion-A2} if
$W \emb \dot B^t_{q,\infty}$ which happens by
\[
     \dot{H}^{-1}_{\sfrac{q}2} \emb \dot{H}^{-1-\sfrac{n}q}_{q} 
\emb \dot{B}^{-1-\sfrac{n}q}_{q,q} \emb \dot{B}^t_{q,\infty}
\]
if $-1-\sfrac{n}q = t$, i.e. if 
$\al + \sfrac{1}p = \frac12 -\frac{n}{2q}$ (see \cite[Theorems
2.7.1 and 5.2.5]{Triebel:FunctionSpaces}).

\smallskip

Finally, using Theorem~\ref{thm:discussion-A1-A3}
\ref{item:discussion-A3j} and Remark~\ref{rem:discussion-A3} the values of 
$\ga$ and $p$ determine $\al$  by 
\begin{equation}\label{eq:alpha-p-n-q-lebesgue}
   \al + \sfrac{1}p = \tfrac12 - \tfrac{n}{2q},
\end{equation}
and then \ref{item:C-zul} and \ref{item:B-zul} are satisfied by the arguments
above.

We sum up the above considerations in the following theorem.
\begin{theorem}\label{thm:appl-lebesgue}
Let $n \geq 2$, $q \in (n,\infty)$, and let $\al\ge0$ and $p\in(2,\infty]$ 
such that \eqref{eq:alpha-p-n-q-lebesgue} holds.
Let $X=\dot{B}^{-1+\sfrac{n}q}_{q,p}(\RR^n)$. Then the Navier-Stokes
equation \eqref{eq:nse} admits a time-local mild solution in
$C([0,\tau), X)$ for every $u_0 \in X^\flat = \overline{\DOMAIN(A)}$ satisfying 
$\nabla\cdot u_0=0$. 
The solution is unique in $C([0,\tau), X) \cap L^p_\al((0,\tau), L^q(\RR^n))$.
If the norm $\norm{u_0}_X$ is sufficiently small, the solution exists
globally. 
\end{theorem}

\begin{remark}
Notice that $X^\flat = X$ in case $p<\infty$ whereas in case $p=\infty$,
$X^\flat$ equals the (homogeneous)
\emph{little Besov space} $\dot{b}^{-1+\sfrac{n}q}_{q,\infty}$
or the (homogeneous) \emph{little Nikolski space} $\dot n^{-1+\sfrac{n}q}_q$ 
(see e.g. \cite{Amann:strong-solvability-NS, Sawada:time-local} for 
the inhomogeneous counterparts).
\end{remark}

\begin{remark}
If we are interested in time--local solutions and use 
Remark~\ref{rem:discussion-finite-time} we are led to $\al+\sfrac{1}p\le 1-\ga$ 
which is equivalent to $2\al+\sfrac{2}p+\sfrac{n}q\le 1$ since 
$\ga=\frac12 +\frac{n}{2q}$. For $\al=0$ we obtain Serrin's uniqueness
 condition $\sfrac{2}p+\sfrac{n}q\le 1$ for weak solutions 
(see e.g. \cite[V.1.5]{Sohr:Buch}). In this context we remark 
that the argument that proved uniqueness in Theorem~\ref{thm:kato-abstrakt}
can be used to show uniqueness of weak solutions 
$u,w\in L^p_\al((0,\tau),L^q(\Om)^n)$ with the same initial value,
but that the assumptions in, e.g., \cite[V.Thm.1.5.1]{Sohr:Buch},
are somewhat weaker and involve energy inequalities.
\end{remark}

\begin{remark}
\begin{enumerate}
\item In case $n=3$, $q>3$, and $p=\infty$ we have 
  $\al=\frac12 - \frac{3}{2q}$ and reobtain a result similar
  to {\sc Cannone} \cite[Theorem 3.3.4]{Cannone:Buch}. 
  There smallness is measured in $X$ but the initial value $u_0$ is taken 
  in $L^3$ and the solution is required to belong to $C_b([0,\tau),L^3)$.
  Since the action $T(t):W=\dot{H}^{-1}_{\sfrac{q}2}\to L^3$ is needed, 
  this leads to the restriction $q<6$ (see \cite{Cannone:Buch}).
\item The general case $n\ge 2$, $q>n$ and $p=\infty$, 
  $\al=\frac12 -\frac{n}{2q}$ is due to {\sc Amann}
  \cite{Amann:strong-solvability-NS} whose proof is similar to taking 
  $W=\dot{H}^{-1-\sfrac{n}q}_q$ for $Z=L^q$. 
  However, \cite{Amann:strong-solvability-NS} also covers the case of
  (sufficiently smooth) domains $\Om\subset\RR^n$, we shall come back to 
  this in Section~\ref{sec:domains-Lq} below.
  Other results on $\RR^n$ with somewhat different approaches are due
  to {\sc Kato} and {\sc Ponce} \cite{KatoPonce}, who used commutator 
  estimates for the bilinear term to achieve existence and regularity
  result for
  initial values in Bessel potential spaces, and to {\sc Koch} and
  {\sc Tataru} \cite{KochTataru}, who proved an existence result for
  initial values
  in $\text{BMO}^{-1}$ and where the structure of proof is more
  involved than via our Theorem~\ref{thm:kato-abstrakt}.   
\item In \cite{Sawada:time-local}, {\sc Sawada} shows existence of 
  time--local mild solutions for divergence--free initial values 
  $u_0\in \overline{\DOMAIN(A)}$ for the inhomogeneous space 
  $X=B^{-1+\sfrac{n}q+\eps}_{q, p}$ where $q\in(n,\infty]$, 
  $p\in[1,\infty]$ and $\eps\in(0,1]$ (thus $X=B^0_{\infty,\infty}$ is
  included). 
  Observe that Theorem~\ref{thm:appl-lebesgue} yields, for $q\in(n,\infty)$ and
  $p\in[\frac{2q}{q-n},\infty]$, local solutions for divergence--free 
  initial values in the space $X=\dot{B}^{-1+\sfrac{n}q}_{q,p}$ 
  (i.e. for $\eps=0$) and global solutions for small initial data 
  (which is not covered by the result in \cite{Sawada:time-local}).
  Moreover, the proof in \cite{Sawada:time-local} relied on a H\"older type
  inequality for products of Besov space functions whereas our proof simply 
  uses the H\"older inequality for the product of two $L^q$--functions. We 
  also remark that we can obtain time-local solutions for 
  the inhomogeneous space $X=B^{-1+\sfrac{n}q}_{q,p}$ by considering
  $\tau<\infty$ and $W=H^{-1}_{q/2}$. We shall discuss the case $q=\infty$ 
  of Sawada's result in Subsection \ref{sec:hoelder} below.

\end{enumerate}
\end{remark}
\bigskip
\subsection{Weak Lebesgue spaces on $\RR^n$}\label{sec:lorentz}
In this section we consider as space $Z$ the weak Lebesgue space 
$L^{q,\infty}$ for a fixed $q\in(n,\infty)$.
For the definition of weak Lebesgue spaces and subsequently used embedding and
interpolation results, see the Appendix in  
Section~\ref{sec:lorentz-def-und-interpol}.
Concerning the Helmholtz projection we remark that, since we are on
$\RR^n$, it commutes with the Laplacian $\Delta$ and that it is
bounded on weak Lebesgue spaces by real interpolation. Similarly,
using \cite[Corollary 6.7.2]{Grafakos} and the interpolation
results of the appendix, the Helmholtz projection is bounded on spaces
$\dot{B}^s_{(q,\infty),p}(\RR^n)$.  
The analysis now follows the lines of Section~\ref{sec:lebesgue}.

For $u,v \in L^{q,\infty}$, clearly $u\otimes v \in L^{\sfrac{q}2,\infty}$ 
and therefore $W := \dot{H}^{-1}_{\sfrac{q}2,\infty}$ guarantees 
$\nabla\cdot (u\otimes v) \in W$. Notice that 
\[
\norm{ T(t) }_{W \to \dot H^\delta_{\sfrac{q}2,\infty}}
\leq c \,t^{-\ga},\quad t>0,
\] 
with $\ga=\tfrac{1+\delta}2$ by bounded analyticity of the semigroup.
By (\ref{eq:homogene-einbettung}) in the proof of Lemma~\ref{lem:claim-two} we have
the embedding $\dot H^{\sfrac{n}q}_{\sfrac{q}2, \infty} \emb L^{q,\infty}$.
Thus $\delta=\sfrac{n}q$ yields the estimate 
$\norm{T(t)}_{W\to Z} \le c\,t^{-\ga}$, $t>0$, required
in Theorem~\ref{thm:discussion-A1-A3}~\ref{item:discussion-A3j} 
and Remark~\ref{rem:discussion-A3} with $\ga=\frac12 +\frac{n}{2q}$. 
Choosing $\al$ and $p$ such that 
\begin{equation}  \label{eq:alpha-p}
  \al + \sfrac1p = 1 -\ga = \tfrac12  - \tfrac{n}{2q},
\end{equation}
condition \ref{item:faltung-bdd} is satisfied. Moreover, letting 
$X :=\dot{B}^{-1+\sfrac{n}q}_{(q,\infty), p}$, the embeddings
\begin{equation}\label{eq:intpol-wbes}
X=(\dot{X}_{1},\dot{X}_{-1})_{\einhalb, p} 
  \quad\text{and}\quad 
(\dot{X}_{1},\dot{X}_{-1})_{\einhalb, \sfrac{p}2} 
  = \dot{B}^{-1+\sfrac{n}q}_{(q,\infty), \sfrac{p}2} \emb X
\end{equation}
hold (for a proof check the corresponding interpolation properties of
vector-valued $\ell^s_p$--spaces \cite[Theorem 1.18.2]{Triebel:interpolation}
and apply a retraction / co-retraction argument). Thus, we can employ 
Theorem~\ref{thm:discussion-A1-A3} in the above setting for verification of the 
assumptions of Theorem~\ref{thm:kato-abstrakt}~\ref{item:C-zul} and 
\ref{item:B-zul}. 

\smallskip

The same arguments that proved \eqref{eq:intpol-wbes} also show
\begin{eqnarray*}
 (X, \dot X_1)_{\al+\sfrac{1}p, 1} = \dot B^s_{(q,\infty), 1} &
  \ \mbox{where}\  & s = -1 + \sfrac{n}q + 2(\al+\sfrac{1}p),\\ 
 (\dot X_{-1}, X)_{2(\al+\sfrac1p), \infty} = \dot B^t_{(q,\infty),\infty}
 & \ \mbox{where}\ & t= 4(\al+\sfrac{1}p)+\sfrac{n}q-3. 
\end{eqnarray*}
By the embedding property \eqref{eq:einschachtelung-schwacher-raeume}
(see Appendix, page~\pageref{eq:einschachtelung-schwacher-raeume})
the first space embeds into $Z = L^{q,\infty}$ for $s = 0$ which
holds by (\ref{eq:alpha-p}). Consequently \ref{item:C-zul} is satisfied.
Similarly, for verification of assumption \ref{item:B-zul},
we observe that $W$ embeds into the second space by Lemma~\ref{lem:claim-two} 
provided that $\al + \sfrac{1}p = \frac12 -\tfrac{n}{2q}$ which again holds 
by (\ref{eq:alpha-p}).
\begin{theorem}\label{thm:appl-lorentz}
Let $n \geq 2$, $q\in(n,\infty)$ and let $\al\ge 0$ and $p\in(2,\infty]$ 
such that (\ref{eq:alpha-p}) holds. 
Let $X := \dot{B}^{-1+\sfrac{n}q}_{(q,\infty), p}$. Then the
Navier-Stokes equation \eqref{eq:nse} admits a time-local mild solution in
$C([0,\tau), X)$ for every $u_0 \in X^\flat = \overline{\DOMAIN(A)}$ satisfying 
$\nabla\cdot u_0 = 0$. The solution is unique in 
$C([0,\tau), X) \cap L^p_\al((0,\tau), L^{q,\infty}(\RR^n))$.
If the norm of $u_0$ is sufficiently small, the solution exists globally.
\end{theorem}
\begin{remark}
Let $q\in(n,\infty)$. In order to give an example of a vector field
that is contained in $\dot{B}^{-1+\sfrac{n}q}_{(q,\infty),\infty}$ but not in 
$B^{-1+\sfrac{n}q}_{q,\infty}$ we use the characterisation of elements that 
are homogeneous of degree $-1$ (see Appendix in
Section~\ref{sec:lorentz-def-und-interpol} where this is shown via 
wavelets). We fix $x_0$ on the unit sphere $S^{n-1}\subset\RR^n$ and
let $w_0(x):=|x-x_0|^{-(n-1)/q}$, $x\in S^{n-1}$. Then 
$w_0\in L^{q,\infty}(S^{n-1})\setminus L^q(S^{n-1})$. We now let 
$v_0:=(1-\Delta)^{(1-\sfrac{n}q)/2}w_0$ where $\Delta$ denotes the 
Laplace-Beltrami operator on $S^{n-1}$. 
By $L^{q,\infty}(S^{n-1})\subseteq
B^0_{(q,\infty),(q,\infty)}(S^{n-1})$ and the lifting 
property we obtain $v_0\in B^{-1+\sfrac{n}q}_{(q,\infty),(q,\infty)}(S^{n-1})$.
We extend $v_0$ by homogeneity of degree $-1$ to the whole of $\RR^n$
and let
\[
u_0(x)=(v_0(x_1,x_2,x_3,\ldots,x_n),-v_0(x_2,x_1,x_3,\ldots,x_n),0,\ldots,0).
\]
Then $\nabla\cdot u_0=0$ and 
$u_0\in\dot{B}^{-1+\sfrac{n}q}_{(q,\infty),\infty}(\RR^n)
 \setminus B^{-1+\sfrac{n}q}_{q,\infty}(\RR^n)$.
\end{remark}
\begin{remark}
In the limit case $q{=}n$ one has $X{=}Z{=}L^{n,\infty}(\RR^n)$ and
$W := \dot{H}^{-1}_{\sfrac{n}2,\infty}$. In this setting, 
existence and uniqueness of 
solutions in $L^\infty((0,\tau),L^{n,\infty}(\RR^n))$ 
has been shown by {\sc Meyer} \cite[Theorem 18.2]{Meyer:Navier-Stokes}.
In our abstract setting we need boundedness of the convolution 
$T(\cdot)\ast:L^\infty(\RR_+,W)\to L^\infty(\RR_+,Z)$, which holds by
Lemma~\ref{lem:L-infty-MR} if $(W,\dot{W}_2)_{\einhalb,\infty}\emb Z$.
By reiteration, the latter condition is equivalent to 
$(\dot{H}^{1-\delta}_{\sfrac{n}2,\infty},\dot{H}^{1+\delta}_{\sfrac{n}2,\infty})_{\einhalb,\infty}
\emb L^{n,\infty}$. This embedding, however, holds by 
$\dot{H}^{1\pm\delta}_{\sfrac{n}2,\infty}\emb L^{n/(1\mp\delta),\infty}$
(see (\ref{eq:homogene-einbettung}) in the proof of Lemma~\ref{lem:claim-two})
and another reiteration identity:
$(L^{n/(1+\delta),\infty},L^{n/(1-\delta),\infty})_{\einhalb,\infty}=L^{n,\infty}$.

\end{remark}

\subsection{Morrey spaces on $\RR^n$}\label{sec:morrey}

In this section we consider as space $Z$ the Morrey space 
$\Morrey^{q,\la}(\RR^n)$ for fixed $q\in(n,\infty)$ and $\la\in(0,\sfrac{n}q)$.
For the definition and some basic properties of Morrey spaces, see the Appendix in  
Section~\ref{sec:lorentz-def-und-interpol}.

\smallskip

For $u,v \in \Morrey^{q,\la}(\RR^n)$, we have 
$u\otimes v \in \Morrey^{\sfrac{q}2,2\la}(\RR^n)$ by H\"older's
inequality.

Therefore, $\nabla{\cdot}(u\otimes v) \in \dot\Morrey^{\sfrac{q}2,2\la,-1}
=: W$ (see the Appendix for the definition of this space).
Observe that, by Calder\'on-Zygmund theory, the Helmholtz projection is
bounded in $W$ and that it commutes with $\Delta$.
Since $W$ equals the homogeneous extrapolation space
$(\Morrey^{\sfrac{q}2,2\la}(\RR^n))_{-\einhalb}^{\displaystyle\cdot}$ with respect to
$-\Delta$ and since $\dot \Morrey^{\sfrac{q}2,2\la,\delta}(\RR^n)$ equals
the homogeneous fractional domain space with respect to $-\Delta$ we have,
by bounded analyticity of the semigroup $T(\cdot)$, 
\[
\norm{ T(t) }_{W \to \dot\Morrey^{\sfrac{q}2,2\la,\delta}(\RR^n)}
\leq c \,t^{-\ga},\quad t>0,
\] 
with $\ga = \frac{1+\delta}2$.

\smallskip

By the properties of the Riesz potential operator (see 
Proposition~\ref{prop:sobolev-emb-morrey} in the Appendix below) we have 
\begin{equation}\label{eq:morrey-emb}
 \dot\Morrey^{\sfrac{q}2,2\la,\delta} \emb \Morrey^{q,\la}
 \qquad\text{for } \delta = \la.
\end{equation}
We conclude $\norm{T(t)}_{W\to Z} \le c\,t^{-\ga}$, $t>0$, where 
$\ga = \frac{1+\la}2$.
Now we use Theorem~\ref{thm:discussion-A1-A3}~\ref{item:discussion-A3j} 
and Remark~\ref{rem:discussion-A3} and choose $\al$ and $p$ such that 
\begin{equation}  \label{eq:alpha-p-morrey}
  \al + \sfrac1p = 1 -\ga = \tfrac{1-\la}{2},
\end{equation}
and condition \ref{item:faltung-bdd} is satisfied. We observe that 
\begin{equation}\label{eq:morrey-X}
 (\dot\Morrey^{q,\la,-2}, \Morrey^{q,\la})_{\ga,p}
  =(\dot\Morrey^{q,\la,-1}, \Morrey^{q,\la})_{ \la , p }
\end{equation}
by reiteration -- using the fact that the dotted spaces are
homogeneous extrapolation spaces for $-\Delta$ in $\Morrey^{q,\la}$.
Denoting the space in (\ref{eq:morrey-X}) by $X$, we clearly have
$X=(\dot{X}_{1},\dot{X}_{-1})_{\einhalb, p}$ by reiteration, hence
also $(\dot{X}_{1},\dot{X}_{-1})_{\einhalb, \sfrac{p}2} \emb X$.
Thus, we can employ Theorem~\ref{thm:discussion-A1-A3} in the above setting
for verification of the assumptions of
Theorem~\ref{thm:kato-abstrakt}~\ref{item:C-zul} and \ref{item:B-zul}. 

\smallskip

We observe that by analyticity of the semigroup and the very
definition of $X$ one has $\norm{T(t)}_{X\to Z} \le c\,t^{-(1-\ga)}$
for $t>0$ and we 
 recall $1-\ga = \al+\sfrac1p$. 
Now Proposition~\ref{prop:simplification-A1-A2} gives $(X,\dot{X}_1)_{\al+\sfrac1p,1}\emb Z$.
Moreover, lifting the embedding \eqref{eq:morrey-emb} by $-(1{+}\la)$
yields
\[
W \emb \dot\Morrey^{q,\la,-(1{+}\la)},
\]
and interpolation yields $\norm{T(t)}_{W\to X}\le c\,t^{-\la}$ for $t>0$.
By Proposition~\ref{prop:simplification-A1-A2} and (\ref{eq:alpha-p-morrey}) we thus have 
$W \emb (\dot{X}_{-1}(-\Delta),X)_{1-\la,\infty}
=(\dot{X}_{-1}(-\Delta),X)_{2(\al+\sfrac1p),\infty}$. 
We have verified the remaining assumptions
\ref{item:C-zul} and \ref{item:B-zul} and obtain the following result.
\begin{theorem}\label{thm:appl-morrey}
Let $n \geq 2$, $\la\in(0,\sfrac{n}q)$, $q\in(n,\infty)$ and let $\al\ge 0$ and 
$p\in(2,\infty]$ such that (\ref{eq:alpha-p-morrey}) holds. 
Let $X$ denote the space in \eqref{eq:morrey-X}. Then the
Navier-Stokes equation \eqref{eq:nse} admits a time-local mild solution in
$C([0,\tau), X)$ for every $u_0 \in X^\flat = \overline{\DOMAIN(A)}$ satisfying 
$\nabla\cdot u_0 = 0$. The solution is unique in 
$C([0,\tau), X) \cap L^p_\al((0,\tau), \Morrey^{q,\la}(\RR^n))$.
If the norm of $u_0$ in $X$ is sufficiently small, the solution exists 
globally.
\end{theorem}

Notice that in the upper limit case $\la=\sfrac{n}q$ in which
$\Morrey^{q,\la}(\RR^n) = L^q(\RR^n)$ condition (\ref{eq:alpha-p-morrey})
becomes exactly condition (\ref{eq:alpha-p-n-q-lebesgue}) we already
found in the case of Lebesgue spaces. The limit case $\la=0$ in which
$\Morrey^{q,\la}(\RR^n) = \text{BMO}$ does not seem to be suited for our 
approach via Theorem~\ref{thm:kato-abstrakt}. In \cite{KochTataru} a somewhat 
different approach allows to take $X=\text{BMO}^{-1}$.  
\begin{remark}
Many authors discussed Navier-Stokes equations in Morrey spaces, see
\cite{Meyer:Navier-Stokes} for an overview. 
The closest result to our theorem is due to {\sc Kozono} and 
{\sc Yamazaki} \cite{KozonoYamazaki} who first introduced real
interpolation spaces of Morrey spaces and of local Morrey spaces (for
the latter, the radii in the definition are restricted to $r\in (0,1]$).
In their notation, our space $X = (\dot\Morrey^{q,\la,-1}, \Morrey^{q,\la})_{ \la , p }$
would be called $\cN^{\la-1}_{\sfrac{n}\la, q, p}(\RR^n)$.
Their main result for Morrey spaces  \cite[Theorem 3]{KozonoYamazaki}
uses $\al=\einhalb-\sfrac\la2$, and $X = \cN^{\la-1}_{\sfrac{n}\la, q, \infty}$ 
which is a special case of our result  for $p=\infty$ but the auxiliary space
there is $\Morrey^{2q, \sfrac{\la}2}$ whereas we have $Z=\Morrey^{q,\la}$.
On the other hand, taking $Z=\Morrey^{2q, \sfrac{\la}2}$ we arrive at the same 
conclusion as Kozono and Yamazaki but for the \emph{larger} space 
$\widetilde{X}=(\dot\Morrey^{2q,\sfrac{\la}2,-1},
\Morrey^{2q,\sfrac{\la}2})_{ \la , \infty}$.
Indeed, by Proposition~\ref{prop:sobolev-emb-morrey} we obtain 
$\norm{T(t)}_{\Morrey^{2q,\la/2}\to\Morrey^{q,\la}}\le c\,t^{-\la/4}$ 
for $t>0$, which implies 
$(\dot\Morrey^{q,\la,-1}, \Morrey^{q,\la})_{ \la , \infty } \emb \widetilde{X}$.
\end{remark}
 
\bigskip
\subsection{H\"older spaces on $\RR^n$}\label{sec:hoelder}

We seek for time-local solutions in this case. In view of
Remark~\ref{rem:discussion-finite-time}, we can assume $A$ to be boundedly
invertible which simplifies the calculation of inter- and
extrapolation spaces.

For fixed $\eps\in(0,1)$, we consider $Z := (C^\eps)^n 
= (B^\eps_{\infty, \infty})^n$. Then for $u,v\in Z$, one has 
$u\otimes v \in (C^\eps)^{n\times n}\subset(\dot{C}^\eps)^{n\times n}$
and thus $\nabla\cdot (u\otimes v)$ 
belongs to the space
\[
  W:=\nabla\cdot(\dot{C}^\eps)^{n\times n} 
   :=\{\nabla\cdot(v_{jk}):(v_{jk})\in (\dot{C}^\eps)^{n\times n} \},
\]
which we equip with the natural quotient--like norm 
\[
\norm{(v_k)}_{\nabla\cdot\dot{C}^\eps}
:= \inf \bigl\{ \norm{ (w_{jk}) }_{\dot{C}^\eps}:\nabla\cdot(w_{jk})=(v_k)\bigr\}.
\]
Observe that, in a canonical way, $\nabla\cdot(\dot{C}^\eps)^{n\times n}$
equals $(\nabla\cdot(\dot{C}^\eps)^n)^n$, and that $W$ is a space of 
distributions although $\dot{C}^\eps$ is not.
Since Riesz transforms are bounded on $\dot{C}^\eps$
(see, e.g. \cite[Corollary 6.7.2]{Grafakos}), 
they are bounded on $W$, and therefore the Helmholtz projection is 
bounded on $W$ (the basic idea is that the origin, in which the symbol 
$\xi_k / |\xi|$ is not differentiable, plays no r\^ole when considering 
\emph{homogeneous} Besov spaces).

We claim that $\norm{T(t)}_{W \to Z} \le C\, \max(1,t^{-\einhalb})$, $t>0$.
Denoting by $(S(\cdot))=(G(\cdot)*)$ the heat semigroup on $\RR^n$, we have
by translation invariance, for $t\in(0,\infty)$,
\[
\norm{S(t)}_{W \to W} \le 1
   \quad\text{and} \quad
\norm{S(t)}_{W\to \dot{C}^\eps} \le C \, t^{-\einhalb}, 
\]
the latter by writing $S(t)\sum_j\partial_jw_j=\sum_j (\partial_jG(t))*w_j$
and using the fact that $\norm{ \partial_jG(t) }_{L^1}\le C\,t^{-\einhalb}$.
Consequently, 
$\norm{T(t)}_{W\to W\cap \dot{C}^\eps} \le C \, \max(t^{-\einhalb},1)$,
and it rests to show $W\cap\dot{C}^\eps\emb Z$, which in turn follows from
$W\cap \dot{C}^\eps\emb L^\infty$. To this end we observe that any $f\in W$ belongs
to $\dot{B}^{\eps-1}_{\infty,\infty}$ and thus has a Littlewood--Paley decomposition 
$f=\sum_{k\in\ZZ} f_k$, for which we obtain
\begin{eqnarray*}
       \Bignorm{ \sum_{k\in\ZZ} f_k }_\infty 
&\le&  \sum_{k\in\ZZ} \bignorm{f_k}_\infty\\ 
& = &  \sum_{k\ge0} 2^{-k \eps } \,\bigl(2^{k \eps }\norm{ f_k }_\infty\bigr)
           + \sum_{k<0} 2^{-k(\eps-1)}\,\bigl(2^{k(\eps-1)}\norm{ f_k }_\infty\bigr)\\
&\le&  \Big(\sum_{k\ge0} 2^{-k\eps}\Big)\, \bignorm{f}_{\dot B^\eps_{\infty,\infty}} 
           + \Big(\sum_{k\ge0} 2^{-k(1-\eps)}\Big)\,
               \bignorm{f}_{\dot B^{\eps-1}_{\infty,\infty}}.
\end{eqnarray*}

\smallskip

Since we are on a finite time interval we can choose 
$\ga =\tfrac{1+\delta}{2} > \einhalb$ where $\delta>0$ is small.
Then $\norm{T(t)}_{W\to Z} \leq c \, t^{-\ga}$ on $(0,\tau)$
for some $c = c_{\ga,\tau} >0$. Such choice of $\ga$ implies 
$\al+\sfrac{1}p=1-\ga=\tfrac{1-\delta}{2}<\einhalb$ whence we obtain for 
$p$ the range $[\tfrac{2}{1-\delta},\infty]$. For such a $p$ we let 
$X=B^{-2(\al+\sfrac{1}p)+\eps}_{\infty, p}=B^{-1+\delta+\eps}_{\infty,p}$.
Then 
\[
({X}_{-1}, {X}_{1})_{\einhalb, p} = X \
   \quad\text{and} \quad
({X}_{-1}, {X}_{1})_{\einhalb, \sfrac{p}2} =
B^{-1+\delta+\eps}_{\infty, \sfrac{p}2} \emb X.
\]
Moreover,
\[
(X,{X}_1)_{\al+\sfrac{1}p, 1} = {B}^{\eps+\delta}_{\infty,1} \emb Z 
   \quad\text{and} \quad
W \emb \dot{B}^{-1+\eps}_{\infty,\infty}\emb 
B^{-1+\eps}_{\infty,\infty}\emb B^{-1-\delta+\eps}_{\infty,\infty} 
\]
where the latter space equals $({X}_{-1}, X)_{2(\al+\sfrac{1}p), \infty}$ 
(recall $2(\al+\sfrac{1}p)=1-\delta$).
We refer to, e.g., \cite[Theorem 2.8.1]{Triebel:interpolation}. 
Therefore, the assumptions of Theorem~\ref{thm:discussion-A1-A3}~\ref{item:discussion-A1} and 
\ref{item:discussion-A2} are satisfied, as well.

\begin{theorem}\label{thm:appl-hoelder} 
Let $n \geq 2$, $p \in (2,\infty]$, and $\eps\in(0,1)$. Let $\al>0$ be such
that $\al+\sfrac{1}p < \einhalb$. Then the Navier-Stokes
equation \eqref{eq:nse} admits a time-local mild solution in
$C([0,\tau), X)$ for every divergence--free 
$u_0 \in X = B^{-2(\al+\sfrac{1}p)+\eps}_{\infty, p}(\RR^n)$, which is unique 
in the space $C([0,\tau), X) \cap L^p_\al((0,\tau), C^\eps(\RR^n))$.
\end{theorem}

\begin{remark}
The result by {\sc Sawada} \cite{Sawada:time-local} also covers time--local solutions for 
initial values in spaces $B^{-1+\eps}_{\infty,p}$ for $p$ up to $\infty$.
However, the space for uniqueness does not involve $L^p_\al(C^\eps)$ but 
$L^\infty_{\beta}$--spaces with values in certain Besov spaces. This is due to the 
fact that the key stone in \cite{Sawada:time-local} is a H\"older type inequality 
for products in (inhomogeneous) Besov spaces which is proved there by means of 
Littlewood-Paley decomposition and paraproducts. Our proof uses the simple product 
inequality in $C^\eps$ instead, and we obtain the second index $p$ in $X$ by taking
$L^p$ in time. So, in our proof, improvement comes from a better understanding of 
the \emph{linear} ingredients for the problem whereas in \cite{Sawada:time-local} it 
comes from a new insight for the \emph{non-linearity}. We remark that 
\cite{Sawada:time-local} includes the case $\eps=1$.
\end{remark}

\bigskip

\subsection{Arbitrary domains in $\RR^3$}\label{sec:domains-L2}

To our knowledge, there are two results in the literature on mild solutions
of the Navier--Stokes equations on arbitrary domains $\Om\subseteq\RR^3$, due
to {\sc Sohr} \cite[Theorem V.4.2.2]{Sohr:Buch} and {\sc Monniaux}
\cite[Theorem 3.5]{Monniaux:arbitrary}. Our results allow to discuss both
approaches, to compare them, and to improve them.

Let $\Om\subseteq\RR^3$ be an arbitrary domain. 
Since there is no regularity assumed for $\partial\Om$, existence of the
Stokes semigroup $(T(t))=(e^{-tA})$ is only guaranteed in $L^2_\si(\Om)$ or 
in interpolation and extrapolation spaces that are associated to 
$L^2_\si(\Om)$ and the Stokes operator $A$.

Since we need the action of $(T(t))$ in $W$ we take $W:=\dot{D}(A^{-\einhalb})
=\dot{\VV}^{-1}_2(\Om)$ (see Section~\ref{sec:prelim}). On $\RR^n$, this would 
correspond to the space $\dot{H}^{-1}_2$, but now we have to pay more 
attention to the Helmholtz projection and $W$ has to be a space of 
divergence--free vectors. We observe that $u,v\in L^4(\Om)^3$ implies 
$u\otimes v\in L^2(\Om)^{3\times 3}$,
$\nabla\cdot(u\otimes v)\in \dot{W}^{-1}_2(\Om)^3$, and finally 
$P\nabla\cdot(u\otimes v)\in\dot{\VV}^{-1}_2(\Om)=W$ by Section~\ref{sec:prelim}.

Since we have Dirichlet boundary conditions, 
$\dot{D}(A^\einhalb)=\dot{\VV}_2\subseteq \dot{W}^1_{2,0}(\Om)^3$ embeds into $L^6(\Om)^3$, 
and by self-adjointness of $A$ and (complex) interpolation we obtain 
$\dot{D}(A^{\sfrac{1}4})\emb L^3(\Om)^3$ and 
$\dot{D}(A^{\sfrac{3}8})\emb L^4(\Om)^3$. Thus, also $u,v\in \dot{D}(A^{\sfrac38})$ implies
$P\nabla\cdot(u\otimes v)\in W$, and $\dot{D}(A^{\sfrac{1}4})$ might be the 
right space of initial values if we seek for global solutions.

\medskip
For $Z\in\{\dot{D}(A^{\sfrac38}),L^4(\Om)^3\}$ we now clearly have  
\[
 \norm{T(t)}_{W\to Z}\le c\norm{T(t)}_{\dot{D}(A^{-\einhalb})\to\dot{D}(A^{\sfrac{3}8})} 
 \leq c \, t^{-\einhalb-\sfrac{3}8},\quad t>0,
\]
i.e. $\gamma=\sfrac{7}8$.

By Theorem~\ref{thm:discussion-A1-A3} we hence should have $\al+\sfrac{1}p=\sfrac{1}8$. 
For inhomogeneities $f_j\in L^{p_j}_{\beta_j}(\RR_+,W^{(j)})$ the condition 
$1+\al+\sfrac{1}p=\gamma_j+\beta_j+\sfrac{1}{p_j}$ then reads 
$\gamma_j+\beta_j+\sfrac{1}{p_j}=\frac{9}8$ where $\gamma_j$ 
is such that  
$\norm{T(t)}_{W^{(j)}\to Z}\le c\,t^{-\gamma_j}$, $t>0$.
Suppose that $X$ is a Banach space satisfying 
\begin{equation}\label{eq:sohrs-X}
 (L^2_\si(\Om),\dot{\VV})_{\einhalb,4}\emb X\emb (L^2_\si(\Om),\dot{\VV})_{\einhalb,8},
\end{equation}
or, for some $p\in(8,\infty]$, 
\begin{equation}\label{eq:allg-sohrs-X}
 (L^2_\si(\Om),\dot{\VV})_{\einhalb,1}\emb X\emb (L^2_\si(\Om),\dot{\VV})_{\einhalb,p},
\end{equation}
and in which the Stokes semigroup acts as a bounded analytic semigroup.
One obtains the pointwise norm estimates 
$\norm{T(t)}_{X\to Z} \le c\, t^{-\sfrac18}$ and
$\norm{T(t)}_{W\to X} \le c\, t^{-\sfrac34}$, see e.g. 
\cite[Lemma 1.12]{HaakKunstmann:weighted-Lp-admiss}.
By reiteration, equation \eqref{eq:sohrs-X} can be reformulated as
\[
 (\dot{X}_{-1},\dot{X}_1)_{\einhalb,4}
\emb X \emb 
 (\dot{X}_{-1},\dot{X}_1)_{\einhalb,8}
\]
and a similar reformulation is possible for \eqref{eq:allg-sohrs-X}.

Thus, we obtain
\begin{theorem}\label{thm:sohr}
Let $Z\in\{\dot{D}(A^{\sfrac38}),L^4(\Om)^3\}$ and $p\in[8,\infty]$. 
Suppose that $X$ is a Banach space satisfying \eqref{eq:sohrs-X}
if 
$p=8$ and \eqref{eq:allg-sohrs-X} 
if 
$p>8$. For any initial value 
$u_0\in X^\flat$ and any $f=f_0+\nabla\cdot F$ with 
$f_0\in L^{p_1}_{\beta_1}(\RR_+,L^2(\Om)^3)$ and 
$F\in L^{p_2}_{\beta_2}(\RR_+,L^2(\Om)^{3\times 3})$, where 
$\beta_j\ge0$, $p_j\in[1,\infty]$
with $\beta_1+\sfrac{1}{p_1}=\sfrac{3}4$
and $\beta_2+\sfrac{1}{p_2}=\sfrac{1}4$, there exists a unique mild solution 
$u$ to the Navier--Stokes equation \eqref{eq:nse} satisfying 
\begin{equation}\label{eq:sohr-prop}
 u\in C\bigl([0,\tau),X\bigr) 
      \; \cap \;   L^p_\al\bigl((0,\tau),Z\bigr)
\end{equation}
where $\al+\sfrac{1}p=\sfrac{1}8$ and $\tau$ depends only on the norms 
$\norm{ u_0 }_X$, $\norm{ f_0 }_{L^{p_1}_{\beta_1}(L^2)}$, 
$\norm{ F }_{L^{p_2}_{\beta_2}(L^2)}$. We 
have $\tau=\infty$ if these norms are sufficiently small.
\end{theorem}
\begin{proof}
We have $f_1=\PP f_0$, $f_2=P\nabla\cdot F$. Taking $W^{(1)}=L^2_\si$ and $W^{(2)}=W$ 
and observing $1+\al+\sfrac{1}p=\sfrac{9}8$,
we only have to check $\norm{T(t)}_{L^2_\si\to Z}\le c\,t^{-\sfrac38}$, i.e. 
$\gamma_1=\sfrac{3}8$, recall $\norm{T(t)}_{W\to Z}\le c\,t^{-\sfrac78}$,
i.e. $\gamma_2=\sfrac{7}8$, and observe 
$\sfrac{9}8-\gamma_1=\sfrac{3}4=\beta_1+\sfrac{1}{p_1}$,
$\sfrac{9}8-\gamma_2=\sfrac{1}4=\beta_2+\sfrac{1}{p_2}$.
\end{proof}
\begin{remark}\label{rem:sohr-max-space}
The following, which takes up an observation from \cite{Cannone:Buch}
shows that, for the choice of $Z=\dot{D}(A^{\sfrac38})$, the space 
$X=(L^2_\si(\Om),\dot{\VV})_{\einhalb,p}$ is maximal for the result:
\begin{eqnarray*}\label{eq:sohr-calc-X}
 \norm{ T(\cdot)x }_{L^p_\al(\dot{D}(A^{\sfrac{3}8}))}
&=&\norm{ t\mapsto t^{\sfrac{1}8}AT(t)x }_{L^p(\RR_+,\frac{dt}{t},\dot{D}(A^{-\sfrac{5}8}))} \\
&\sim& \norm{ x }_{(\dot{D}(A^{-\sfrac{5}8}),\dot{D}(A^{\sfrac{3}8}))_{\sfrac{7}8,p}}
\sim \norm{ x }_{(L^2_\si(\Om),\dot{\VV})_{\einhalb,p}}.
\end{eqnarray*}

\end{remark}
\begin{remark}
Sohr's result (\cite[Theorem V.4.2.2]{Sohr:Buch}) has 
$Z=L^4(\Om)^3$, $p=8$, $\al=0$, $\beta_1=\beta_2=0$ and $p_1=\sfrac{4}3$,
$p_2=4$. It is remarkable that \eqref{eq:sohrs-X} does not allow to take 
$X=\dot{D}(A^{\sfrac{1}4})$. In fact, Sohr takes \emph{weak} solutions
$u$ of (\ref{eq:nse}) which always satisfy 
$u\in L^\infty_{loc}([0,\eta),L^2_\si(\Om))\cap 
L^2_{loc}([0,\eta),\VV_2(\Om))$.
Observe that the space $X=\dot{D}(A^{\sfrac{1}4})$ becomes admissible 
if we choose $p>8$.
\end{remark}
\begin{remark}\label{rem:new-monniaux}
Taking $Z=\dot{D}(A^{\sfrac38})$ and $p=\infty$, $\al=\sfrac{1}8$ in 
Theorem~\ref{thm:sohr} we obtain an improvement of {\sc Monniaux}'s result 
(\cite[Theorem 3.5]{Monniaux:arbitrary}, see also the discussion below). 
Here we may choose $X=(L^2_\si(\Om),\dot{\VV})_{\einhalb,s}$ for $s\in[1,\infty]$.
Thus the maximal space for initial values is 
$X=(L^2_\si(\Om),\dot{\VV})_{\einhalb,\infty}$.
Observe that $X=\dot{D}(A^{\sfrac{1}4})$ for $s=2$ since $A$ is selfadjoint.\\

In \cite{Monniaux:arbitrary} the right hand side is $f=0$. Moreover,
the assertion there only covers time--local solutions. Actually, the space 
$\mathscr{E}_T$ in \cite{Monniaux:arbitrary} is not a Banach space, in general.
Since only time--local solutions are considered in \cite{Monniaux:arbitrary},
the proof can be corrected by replacing $A$ in the definition of the 
norm of $\mathscr{E}_T$ with $\delta+A$ in case $0\in\si(A)$. In this context,
we remark that $\dot{\VV}=\VV$, $W=\VV'$ and $\dot{D}(A^r)=D(A^r)$
for $r>0$ if $0\in\rho(A)$ which happens, e.g., if $\Om$ is bounded.
\end{remark}
We want to discuss the result in \cite{Monniaux:arbitrary} a bit further. The approach 
there corresponds to taking $Z=\VV=D(A^\einhalb)$. For $u,v\in Z$, one has 
$u\cdot\nabla v\in L^{\sfrac{3}2}(\Om)^3$. Dualising 
$D(A^{\sfrac{1}4})\emb L^3(\Om)^3$ yields
$P:L^{\sfrac{3}2}(\Om)^3\to (D(A^{\sfrac{1}4}))'$, and the latter space equals 
$W:=(L^2_\si(\Om),\norm{ (\delta+A)^{-\sfrac{1}4}\cdot })^\sim$ (the 
embedding $\dot{D}(A^{\sfrac{3}8})\emb L^3(\Om)^3$ might be used as well;
then $P:L^{\sfrac{3}2}(\Om)^3\to\dot{D}(A^{-\sfrac{1}4})$, and one could choose 
$W:=\dot{D}(A^{-\sfrac{1}4})$, but the other choice is closer to what is 
actually happening in \cite{Monniaux:arbitrary}). Now clearly 
$\norm{ T(t) }_{W\to Z}\le c \,t^{-\sfrac{3}4}$ on bounded time--intervals. 
In order to satisfy $\al+\sfrac{1}p\le 1-\sfrac{3}4=\sfrac{1}4$ choose 
$p=\infty$ and $\alpha=\sfrac{1}4$. If $X$ is a Banach space satisfying
\[
 (L^2_\si(\Om),\VV)_{\einhalb,1}\emb
 X\emb (L^2_\si(\Om),\VV)_{\einhalb,\infty}
\]
in which the Stokes semigroup acts as a bounded analytic semigroup, then we obtain
\begin{theorem}\label{thm:monniaux}
For any initial value $u_0\in X^\flat$ 
and $f=f_0+\nabla\cdot F$ with $f_0\in L^{p_1}_{\beta_1}(\RR_+,L^2(\Om)^3)$ and 
$F\in L^p_{\al}(\RR_+,L^2(\Om)^{3\times 3})$, where $\al,\beta_1\ge0$,
$p_1\le\infty$, $p<\infty$ and $\al+\sfrac{1}p=\sfrac{1}4$, 
$\beta_1+\sfrac{1}{p_1}=\sfrac{3}4$,
there is a unique mild solution $u$ to the Navier--Stokes equation (\ref{eq:nse}) 
satisfying 
\begin{equation}\label{eq:monniaux-prop}
 u\in C_b([0,\tau),X)
      \cap L^p_{\al}((0,\tau),\VV)
\end{equation}
where $\tau$ depends only on the norms $\norm{ u_0 }_X$, 
$\norm{ f_0 }_{L^{\beta_1}_{\beta_1}(L^2)}$, $\norm{ F }_{L^p_\al(L^2)}$. 
\end{theorem}
\begin{proof}
We have $f_1=\PP f_0$, $f_2=P\nabla\cdot F$, $W^{(1)}=L^2_\si$, $W^{(2)}=\VV^{-1}_2$.
Observing $\|T(t)\|_{L^2_\si\to\VV}\le ct^{-\einhalb}$, 
$\|T(t)\|_{\VV^{-1}_2\to \VV}\le ct^{-1}$ on finite time intervals, this leads 
to $\ga_1=\einhalb$ and $\ga_2=1$. Notice that 
$1+\al+\sfrac{1}p-\einhalb=\sfrac{5}4-\einhalb=\sfrac{3}4=\beta_1+\sfrac{1}{p_1}$.
However, now that $\ga_2=1$ we need maximal $L^p_\al$-regularity (on finite time 
intervals, see discussion in Section~\ref{sec:abstract-kato}) for the Stokes operator
$A$ in $\VV^{-1}_2$, which holds since $\VV^{-1}_2$
is a Hilbert space, $p<\infty$, and $0\le\al\le\sfrac{1}4-\sfrac{1}p$.
\end{proof}
\begin{remark}
As mentioned before, \cite{Monniaux:arbitrary} has $f=0$. Observe that, although the
space $Z$ is different, the conditions on the right hand side $f$ are the same as in 
Theorem~\ref{thm:sohr}, but that $\ga_2=1$ led to the restriction $\beta_2=\al$, $p_2=p<\infty$,
since we need the continuous action 
$T(\cdot)*:L^{p_2}_{\beta_2}(\VV^{-1}_2)\to L^p_\al(\VV_2)$. 
If $F=0$ we can admit $p=\infty$ in the assertion.
Observe also that it was essential for the argument to use the \emph{inhomogeneous} 
space $\VV$, which in turn restricts the result to time-local solutions.
\end{remark}

\subsection{Domains which admit an $L^q$-theory}\label{sec:domains-Lq}
In this subsection $\Om\subseteq\RR^n$ is a domain 
for which we assume additionally that,
for some $q_0\in(2,\infty)$, the Helmholtz projection is bounded in 
$L^{q_0}(\Om)^n$ and the Stokes semigroup is bounded analytic 
in $L^{q_0}_\si(\Om)$ (see the end of Section~\ref{sec:prelim}). 
We distinguish two cases.

\smallskip

\textbf{Case I, $n=3$ and $q_0\in(2,4]$:} 
We start with a preparation. By interpolating the semigroup estimates
$\norm{ T(t) }\le c \,t^{-\einhalb}$ for the action $T(t):L^2_\si\to L^6$ and 
$\norm{ T(t) }\le c$ for the action $T(t):L^{q_0}_\si\to L^{q_0}$ one obtains
\begin{equation}\label{eq:Tt-Lq-L4}
 \norm{ T(t) }_{L^q_\si\to L^4}\le c \,t^{-\sfrac{3}{2}(\sfrac{1}q-\sfrac{1}4)}, \quad t>0,
\end{equation}
where $\theta$ is determined by $\frac{1}{4}=\frac{1-\theta}{6}+\frac{\theta}{q_0}$ and 
$q$ satisfies $\frac{1}{q}=\frac{1-\theta}{2}+\frac{\theta}{q_0}$ (observe that the
negative $t$-exponent equals $\delta:=\frac{3}{2}(\frac{1}{q}-\frac{1}{4})
                                     =\frac{3}{2}\frac{1-\theta}{3}
                                     =\frac{1-\theta}{2}+\theta\cdot 0$).

Choose $W:=\dot{\VV}^{-1}_2(\Om)$ and $Z:=L^4(\Om)^3$ as in the previous 
subsection for Theorem~\ref{thm:sohr}. Then still $\al+\sfrac{1}p=\sfrac{1}8$. 
We want to find spaces $X$ associated to $L^q_\si(\Om)$ and $A_q$.
To this end we let $\widetilde Z:=(L^q_\si(\Om),\dot{D}(A_q))_{\delta,1}$.
By \eqref{eq:Tt-Lq-L4} and Proposition~\ref{prop:simplification-A1-A2} we know that 
$\widetilde Z\emb L^4(\Om)^3=Z$. Using the Stokes semigroup in $L^q_\si(\Om)$
we calculate (as in Remark~\ref{rem:sohr-max-space})
\[
 \norm{ T(\cdot)x }_{L^p_\al(\widetilde Z)}
 =\norm{ t\mapsto t^{\sfrac{1}8}AT(t)x }_{L^p(\frac{dt}{t},(\widetilde Z)^{\displaystyle\cdot{}}_{-1})}
 \sim\norm{ x }_{((\widetilde Z)^{\displaystyle\cdot{}}_{-1},\widetilde Z)_{\sfrac{7}8,p}}.
\]
Clearly, $(\widetilde Z)^{\displaystyle\cdot{}}_{-1}=(\dot{D}(A_q^{-1}),L^q_\si(\Om))_{\delta,1}$, 
and by reiteration,
\[
  ((\widetilde Z)^{\displaystyle\cdot{}}_{-1},\widetilde Z)_{\sfrac{7}8,p} 
 =(\dot{D}(A_q^{-\sfrac{1}2}),
   \dot{D}(A_q^{\sfrac{1}2}))_{\frac{3}{2q},p}.
\]
Indeed, observe that 
$\frac{1}{8}(\delta-1)+\frac{7}{8}\delta=\delta-\frac{1}{8}=\frac{3}{2q}-\frac{1}{2}$
and $(1-\frac{3}{2q})(-\frac{1}{2})+\frac{3}{2q}\cdot\frac{1}{2}=-\frac{1}{2}+\frac{3}{2q}$.

As an illustration we remark that, for $\Om=\RR^3$, this space equals
the homogeneous divergence--free Besov space 
$\dot{B}^{-1+\frac{3}{q}}_{q,p,\si}$.
\begin{theorem}\label{thm:Lq-domain-I}
Suppose $q_0$ and $q$ are as above. Let $\al\ge0$, $p\in[8,\infty]$ such that
$\al+\sfrac{1}p=\sfrac{1}8$, and let 
$X:=(\dot{D}(A_q^{-\sfrac{1}2}),\dot{D}(A_q^{\sfrac{1}2}))_{\frac{3}{2q},p}$.
For any initial value $u_0\in X^\flat$ and any $f=f_0+\nabla\cdot F$ with 
$f_0\in L^{p_1}_{\beta_1}(\RR_+,L^2(\Om)^3)$ and 
$F\in L^{p_2}_{\beta_2}(\RR_+,L^2(\Om)^{3\times 3})$, where 
$\beta_j\ge0$, $p_j\in[1,\infty]$
with $\beta_1+\sfrac{1}{p_1}=\sfrac{3}4$
and $\beta_2+\sfrac{1}{p_2}=\sfrac{1}4$,
there is a unique mild solution to the Navier-Stokes equation (\ref{eq:nse}) satisfying
\[
 u\in C_b([0,\tau),X)\cap L^p_\al((0,\tau),L^4(\Om)^3)
\]
where $\tau$ depends only on the norms $\norm{u_0}_X$, 
$\norm{f_0}_{L^{p_1}_{\beta_1}(L^2)}$, $\norm{F}_{L^{p_2}_{\beta_2}(L^2)}$, 
and we have $\tau=\infty$ if these norms are sufficiently small.
\end{theorem}
\begin{remark}
Concerning the relation of $q$ and $q_0$ we remark that a calculation shows 
$\frac{1}{q}=\frac{1}{4}(2-\frac{q_0-2}{6-q_0})=\frac{1}{4}\frac{12-3q_0+2}{6-q_0}$.
Hence we have $q=q_0$ for $q_0\in\{2,4\}$, and the special cases 
$q=\sfrac{12}5$ for $q_0=3$ and $q=3$ for $q_0=\sfrac{18}5$. 
We did not use any further properties besides boundedness
of the Helmholtz projection in $L^{q_0}$ and bounded analyticity of the Stokes semigroup
in $L^{q_0}_\si$. Once one has
\begin{equation}\label{eq:Tt-Lq1-L4}
 \norm{T(t)}_{L^{q_1}_\si\to L^4}\le c \,t^{-\sfrac{3}{2}(\sfrac{1}{q_1}-\sfrac{1}4)}, 
 \quad t>0,
\end{equation}
for some $q_1\in(q,q_0]$, in Theorem~\ref{thm:Lq-domain-I} the spaces 
$(\dot{D}(A_{q_1}^{-\sfrac{1}2}), \dot{D}(A_{q_1}^{\sfrac{1}2}))_{\frac{3}{2q_1},p}$
can be used in place of the spaces $(\dot{D}(A_q^{-\sfrac{1}2}),
\dot{D}(A_q^{\sfrac{1}2}))_{\frac{3}{2q},p}$.
In this context, we remark that the spaces 
$(\dot{D}(A_q^{-\sfrac{1}2}),\dot{D}(A_q^{\sfrac{1}2}))_{\frac{3}{2q},p}$
actually grow with $q$. Since $A^{-\einhalb}$ is an isometry, it is sufficient to show 
that $(L^q_\si(\Om),\dot{D}(A_q))_{\frac{3}{2q},p}$ grows with $q$. To see this,
we let $q_1\in(q,q_0]$, write out the norms and use the semigroup property 
\begin{eqnarray*}
 \norm{x}_{(L^{q_1}_\si(\Om),\dot{D}(A_{q_1}))_{\frac{3}{2q_1},p}}
 &\sim & \norm{t\mapsto t^{1-\frac{3}{2q_1}}AT(t)x}_{L^p(\RR_+,\frac{dt}{t},L^{q_1})}\\
 &\sim &\norm{t\mapsto t^{1-\frac{3}{2q_1}}T(t)AT(t)x}_{L^p(\RR_+,\frac{dt}{t},L^{q_1})}\\
 &\le & c \norm{t\mapsto t^{1-\frac{3}{2q}}AT(t)x}_{L^p(\RR_+,\frac{dt}{t},L^{q})},
\end{eqnarray*}
where we used
\begin{equation}\label{eq:q-q1-est}
 \norm{T(t)}_{L^q_\si\to L^{q_1}}\le c t^{-\frac{3}{2}(\frac{1}{q}-\frac{1}{q_1})},
 \quad t>0,
\end{equation}
in the last step, which in turn follows by interpolation of the action
$T(t):L^{q_1}_\si\to L^{q_1}$ and $T(t):L^2_\si\to L^{q_1}$, recall
that $q_1\le 4<6$.
It is clear that, besides bounded analytic action of the Stokes semigroup in 
$L^q_\si(\Om)$ and $L^{q_1}_\si(\Om)$, the estimate \eqref{eq:q-q1-est} is all
that is needed to prove the desired inclusion in more general cases.
\end{remark}

\smallskip

\textbf{Case II, $q_0>\max\{n,4\}$:} We let $q:=q_0$ for simplicity of notation.
One can choose $Z:=L^q(\Om)^n$. For $u,v\in Z$ then
$u\otimes v\in L^{\sfrac{q}2}(\Om)^{n\times n}$, 
$\nabla\cdot(u\otimes v)\in \dot{W}^{-1}_{\sfrac{q}2}(\Om)^n$, 
and $P\nabla\cdot(u\otimes v)\in \dot{\VV}^{-1}_{\sfrac{q}2}(\Om)=:W$.
Notice that $\sfrac{q}2>2$ by $q>4$. We now aim at
\begin{equation}\label{eq:cond-q-Tt}
 \norm{ T(t) }_{\dot{\VV}^{-1}_{\sfrac{q}2}(\Om)\to L^q(\Om)^n}
 \le c \,t^{-\frac{1}2-\frac{n}{2q}},\quad t>0,
\end{equation}
i.e. $\ga=\tfrac12-\tfrac{n}{2q}$. Here, $L^q(\Om)^n$ can be replaced by
$L^q_\si(\Om)$, in which space we can use the Stokes semigroup to obtain 
the equivalent condition 
\begin{equation}\label{eq:cond-q-inf}
 \dot{\VV}^{-1}_{\sfrac{q}2}(\Om)
 \emb(\dot{D}(A_q^{-1}), L^q_\si(\Om))_{\frac{1}2-\frac{n}{2q},\infty}
\end{equation}
by Proposition~\ref{prop:simplification-A1-A2}. Dualising 
\eqref{eq:cond-q-Tt} (with $L^q_\si(\Om)$ in place of $L^q(\Om)^n$)
yields as another equivalent condition 
\begin{equation}\label{eq:cond-q'-Tt}
 \norm{ T(t) }_{L^{q'}_\si(\Om)\to \dot{\VV}_{(\sfrac{q}2)'}(\Om)}
 \le c \,t^{-\frac{1}2-\frac{n}{2q}},\quad t>0,
\end{equation}
which in turn can be reformulated by Proposition~\ref{prop:simplification-A1-A2} as
\begin{equation}\label{eq:cond-q'-1}
  (L^{q'}_\si(\Om),\dot{D}(A_{q'}))_{\frac{1}2+\frac{n}{2q},1}
  \emb \dot{\VV}_{(\sfrac{q}2)'}(\Om),
\end{equation}
where we used the Stokes semigroup in $L^{q'}_\si(\Om)$. Another 
reformulation of \eqref{eq:cond-q'-Tt} is the following gradient estimate
\begin{equation}\label{eq:cond-grad}
  \norm{ \nabla T(t)f }_{(\sfrac{q}2)'}
  \le c \,t^{-\frac{1}2-\frac{n}{2q}}\norm{ f }_{q'},
  \quad t>0, f\in L^{q'}_\si(\Om).
\end{equation}
We thus obtain the following new result.
\begin{theorem}\label{thm:Lq-domain-II}
Suppose that $q=q_0$ is as above and assume that one of the equivalent conditions
\eqref{eq:cond-q-Tt}, \eqref{eq:cond-q-inf}, \eqref{eq:cond-q'-Tt}, \eqref{eq:cond-q'-1},
\eqref{eq:cond-grad} holds.
Let $\al\ge0$, $p\in(2,\infty]$ such that
$\al+\sfrac{1}p=\frac{1}{2}-\frac{n}{2q}$, and let 
$X:=(\dot{D}(A_q^{-\sfrac{1}2}),\dot{D}(A_q^{\sfrac{1}2}))_{\frac{n}{2q},p}$.
For any initial value $u_0\in X^\flat$ and $f=f_0+\nabla\cdot F$ with 
$f_0\in L^{p_1}_{\beta_1}(\RR_+,L^2(\Om)^3)$ and 
$F\in L^{p_2}_{\beta_2}(\RR_+,L^2(\Om)^{3\times 3})$, where 
$\beta_j\ge0$, $p_j\in[1,\infty]$
with $\beta_1+\sfrac{1}{p_1}=\sfrac{3}2-\sfrac{n}q$
and $\beta_2+\sfrac{1}{p_2}=1-\sfrac{n}q$
there is a unique mild solution to the Navier-Stokes equation (\ref{eq:nse}) 
satisfying
\[
 u\in C_b([0,\tau),X)\cap L^p_\al((0,\tau),L^q(\Om)^n)
\]
where $\tau$ depends only on the norms $\norm{u_0}_X$, 
$\norm{f_0}_{L^{p_1}_{\beta_1}(L^2)}$, $\norm{F}_{L^{p_2}_{\beta_2}(L^2)}$, 
and we have $\tau=\infty$ if these norms are sufficiently small.
\end{theorem}
\begin{proof}
As mentioned above we have $Z=L^q(\Om)^n$ and $W =\dot{\VV}^{-1}_{\sfrac{q}2}(\Om)$. 
We have $\ga=\frac{1}{2}+\frac{n}{2q}$ which explains the 
condition on $\al+\sfrac{1}p$. Notice that
$W^{(1)} = L^{\sfrac{q}2}_\si(\Om)$ yields $\ga_1 = \tfrac{n}{2q}$
and since $W^{(2)}=W$,  $\ga_2=\ga$. 
To verify \ref{item:A3-inhomogen}, by \ref{thm:discussion-A1-A3} requires
$\beta_1+\sfrac1{p_1} +\ga_1 = \al+\sfrac1p+1$ and
$\al+\sfrac1p+\ga=1$ -- both are guaranteed by the assumptions on 
$\beta_1$/$\beta_2$ 
and $p_1$/$p_2$. 
Moreover, we have 
$(\dot{X}_{-1},\dot{X}_1)_{\einhalb,\sfrac{p}2}\emb(\dot{X}_{-1},\dot{X}_1)_{\einhalb,p}=X$
by reiteration, and an argument as in Remark~\ref{rem:sohr-max-space} shows that
$X$ satisfies \ref{item:C-zul}. Finally, assumption \ref{item:B-zul} 
follows from (\ref{eq:cond-q-inf}) by reiteration. 
\end{proof}

\begin{remark}
If $\Om$ is bounded and $\partial\Om$ is of class $C^{1,1}$ then
Theorem~\ref{thm:Lq-domain-II} may be applied to any $q\in(n,\infty)$. It is
well-known that the Stokes semigroup is bounded analytic in all
$L^q_\si$, $q\in(1,\infty)$.  Moreover, condition \eqref{eq:cond-grad}
is satisfied for any $q\in(n,\infty)$.  This follows from 
\[
 \norm{T(t)}_{L^{q'}_\si\to L^{(\sfrac{q}2)'}_\si}\le c t^{-\frac{n}{2q}},\quad t>0,
\]
and
\begin{equation}\label{eq:grad-q2'}
 \norm{\nabla T(t)}_{L^{(\sfrac{q}2)'}_\si\to L^{(\sfrac{q}2)'}}\le c t^{-\einhalb}, \quad t>0,
\end{equation}
where the latter is due to the fact that $A_{(\sfrac{q}2)'}$ has a
bounded $H^\infty$-calculus in $L^{(\sfrac{q}2)'}_\si(\Om)$
(\cite[Thm.9.17]{KaltonKunstmannWeis}), hence has bounded imaginary
powers,  which leads to
$\dot{D}(A_{(\sfrac{q}2)'}^\einhalb)=[L^{(\sfrac{q}2)'}_\si(\Om),\dot{D}(A_{(\sfrac{q}2)'}]_\einhalb=
\dot{\VV}_{(\sfrac{q}2)'}$ (observe that $0\in\rho(A_{(\sfrac{q}2)'})$
and thus homogeneous and inhomogeneous domain spaces coincide). 
The result on mild solutions in this case is due to 
{\sc Amann} \cite{Amann:strong-solvability-NS}.

An approach to unbounded domains $\Om$ of uniform $C^{1,1}$-type based on
$\widetilde{L}^q(\Om)$-spaces is due to {\sc Farwig}, {\sc Kozono},
and {\sc Sohr} \cite{FarwigKozonoSohr}. Here, ${\widetilde L}^q=L^q+L^2$
for $q\le2$ and
$\widetilde{L}^q=L^q\cap L^2$ for $q>2$. Resorting to these spaces,
difficulties that arise for an $L^q$-theory from the behaviour of $\Om$
at infinity could be overcome. We refer to \cite{FarwigKozonoSohr} for
details.
\end{remark}

\section{{Proof of Theorem~\ref{thm:discussion-A1-A3}}}
\label{sec:proof-of-kato-discussion}

For a sectorial operator in a Banach space $X$, real interpolation spaces 
between $X$ and inhomogeneous spaces $X_k$ are well-studied,
see \cite{BerghLoefstroem,Lunardi:Buch,Triebel:interpolation}.
In this section we provide the results on real interpolation of 
\emph{homogeneous} spaces needed for the proof of
Theorem~\ref{thm:discussion-A1-A3}. The following result is an 
analogue of \cite[Theorem 1.14.2]{Triebel:interpolation}.
\begin{proposition}\label{prop:homogene-interpol-R}
  Let $X$ be a Banach space and $A$ be an injective, sectorial operator
on $X$. Then, for $m\in \NN$, $(X, \dot{X}_m)$ is a quasi-linearisable
interpolation couple in the sense of 
\cite[Definition 1.8.3]{Triebel:interpolation}. Moreover, for $p\in
[1,\infty]$ and $\theta\in(0,1)$, an equivalent norm on 
$(X, \dot{X}_m)_{\theta,p}$ is given by 
$\norm{\la \mapsto \la^{\theta m} A^m (\la+A)^{-m} x}_{L^p(\RR_+, d\la/\la, X)}$.
\end{proposition}
\begin{proof}
Let $E_m := X + \dot{X}_m$. Clearly, $E_m = (I+\AA)^m (\dot{X}_m)$ and 
$\norm{e}_{E_m} = \norm{ \AA^m(I+\AA)^{-m} e}_X$. We borrow a
decomposition technique inspired by 
\cite[Proposition 15.26]{KunstmannWeis:Levico}:
let $a_j$ be defined by $\sum_{j=1}^{2m-1} a_j z^j = (1+z)^{2m} -
(1+z^m)(1+z)^m$. Therefore, setting $z=\la^{\sfrac{1}m} A$, we obtain
for all $x\in E_m$, 
\begin{eqnarray*}
x &=& \biggl[ \la A^m (1+\la^{\sfrac{1}m} A)^{-m} x + \sum_{j=m}^{2m-1} a_j
      \la^\sfrac{j}m A^j (1+\la^{\sfrac{1}m} A)^{-2m} x \biggr]\\
&& + \biggl[ (1+\la^\sfrac{1}mA)^{-m} x + \sum_{j=1}^{m-1} a_j
      \la^\sfrac{j}m A^j (1+\la^{\sfrac{1}m} A)^{-2m} x \biggr]
\end{eqnarray*}
(all operators are bounded). Call the first term in brackets $V_0(\la) x$
and the second one $V_1(\la) x$. A direct calculation shows 
quasi-linearisability. Hence,
\[
K(\la,x,X,\dot{X}_m) \sim \norm{V_0(\la)x}_X +
\la\norm{V_1(\la)x}_{\dot{X}_m}.
\]
Notice that
\[ 
  \la^{-\theta} (\la^\sfrac{1}m A)^m (I+\la^\sfrac{1}m A)^{-m}x
= \la^{-\theta} A^m (\la^{-\sfrac{1}m} + A)^{-m}x
= \tau^{\theta m} A^m (\tau+A)^{-m}x,
\]
by letting $\tau = \la^{-\sfrac{1}m}$. It remains to show
\[
\bignorm{V_0(\la)x}_X + \la \bignorm{V_1(\la)x}_{\dot{X}_m} 
\quad \sim \quad
\bignorm{(\la^\sfrac{1}m A)^m (I+\la^\sfrac{1}m A)^{-m}x}.
\]
For the estimate ``$\leq$'', notice that 
\[
\la^{\tfrac{j+m}m} A^{m+j} (I+\la^\sfrac{1}m A)^{-2m}
=
(\la^\sfrac{1}m A)^j (I+\la^\sfrac{1}m A)^{-m} \cdot \la A^m (I+\la^\sfrac{1}m A)^{-m},
\]
where the first expression is bounded by sectoriality of $A$.

Finally, let $f(\la) := \la A^m(1+\la^\sfrac{1}mA)^{-m}$. Then for $x \in
E_m$ and $x=y+z$ with $y\in X$ and $z \in \dot{X}_m$,
\[
f(\la) x = f(\la) y + f(\la) z = f(\la) y + \la \bigl(\la^{-1}A^{-m}\bigr)
f(\la) A^m z,
\]
whence by sectoriality of $A$,
\[
\norm{f(\la)x}_X 
\leq \norm{f(\la)y}_X + M \la \norm{z}_{\dot{X}_m}.
\]
Taking the infimum over all such decomposition yields
$\norm{f(\la)x}_X  \leq c \, K(\la,x,X,\dot{X}_m)$ and the proof is finished.
\end{proof}

The following result corresponds to 
\cite[Theorem 1.14.5]{Triebel:interpolation} and gives another
equivalent norm on $(X, \dot{X}_m)_{\theta,p}$ in the case
of analytic semigroups. 
We omit the proof since it is identical to the non--homogeneous case. 
\begin{proposition}\label{prop:homogene-interpol-T}
Let $T(\cdot)$ be a bounded and analytic semigroup on a Banach space
$X$ and $-A$ its generator. If $A$ is injective, then for 
$p\in [1,\infty]$ and $\theta\in(0,1)$, an equivalent norm on 
$(X, \dot{X}_m)_{\theta,p}$ is given by 
$\norm{ t^{m(1-\theta)} A^m T(t) x}_{L^p(\RR_+, dt/t, X)}$.
\end{proposition}

The following result is an analogue of 
\cite[Theorem 1.14.3 (a)]{Triebel:interpolation}. 
\begin{lemma}\label{lem:einbettung}
Let $X$ be a Banach space and $A$ be an injective, sectorial operator
on $X$. Then, for $k,j,m\in\ZZ$ with $k< j < m$, we have
\[
     \bigl(\dot{X}_k, \dot{X}_m\bigr)_{\frac{j-k}{m-k},1} 
\emb \dot{X}_j 
\emb \bigl(\dot{X}_k, \dot{X}_m\bigr)_{\frac{j-k}{m-k},\infty}.
\]  
\end{lemma}
\begin{proof}
We can assume $k=0$.
Let $x\in \DOMAIN(A^m) = X \cap \dot{X}_m \subseteq (X,
\dot{X}_m)_{\sfrac{j}m,1}$. Then, by
\cite[1.14.2/(1)]{Triebel:interpolation}, for some constant $c_m$
\[
x = c_m \int_0^\infty (tA)^m (t+A)^{-2m} x\,\tfrac{dt}t.
\]
Therefore, by sectoriality of $A$,
\begin{eqnarray*}
 \bignorm{ A^j x}_X 
&\leq& c_m \int_0^\infty t^{m-1} \bignorm{A^{j+m}(t+A)^{-2m} x}_X \,dt\\
&=& c_m \int_0^\infty  \bignorm{t^{m-j} A^j (t+A)^{-m} \; t^j
    \bigl[A(t+A)^{-1}\bigr]^m  x} \,\tfrac{dt}t\\
&\leq& M c_m \int_0^\infty  \bignorm{t^j \bigl[A(t+A)^{-1}\bigr]^m  x} \,\tfrac{dt}t
\leq  \widetilde{c}_m \bignorm{x}_{(X, \dot{X}_m)_{\sfrac{j}m,1}}.
\end{eqnarray*}
The second embedding follows also by sectoriality of $A$ from
\[
\bignorm{ t^j [A(t+A)^{-1}]^m x}_{X} = 
\bignorm{ t^j A^{m-j} (t+A)^{-m} A^j x}_{X}
\leq \widetilde{M} \norm{x}_{\dot{X}_{j}},
\]
which is true for all $x\in \dot{X}_j$.
\end{proof}

The assertion does hold for arbitrary interpolation indices
$\theta\in(0,1)$ but we shall not introduce fractional
homogeneous spaces (see \cite{HaakHaaseKunstmann,Haase:Buch}) 
since the above version is sufficient for our purposes.

\subsection{Results on assumption {\ref{item:C-zul}}} \label{sec:on-A1}

In this section we discuss boundedness of the map
\begin{equation}\label{eq:C-zul}
  \Psi_\infty: X \to L^p_\al(\RR_+, Z), 
\qquad \Psi_\infty(x) = C   T(\cdot) x 
\end{equation}
for $p \in [1, \infty]$ and $\al\in(-\sfrac{1}p, 1-\sfrac{1}p)$. We
start our considerations with a simple necessary condition for
boundedness of 
$\Psi_\infty$:
 boundedness of the set
\begin{equation} \label{eq:necess-C-zul}
  \{ \la^{1-\al-\sfrac{1}p} C(\lambda+A)^{-1} \suchthat \la >0 \}
\end{equation}
in $B(X, Z)$. Indeed, writing the resolvent of $A$ as Laplace
transform of the semigroup and using H\"older's inequality we have for
$x\in X_1$ and $\la>0$
\begin{eqnarray*}
       \bignorm{ C(\lambda+A)^{-1} x } 
&\leq & \int_{0}^{\infty} \bignorm{t^{\alpha} CT(t) } t^{-\alpha} e^{-\la t} \,dt\\
&\leq & \biggnorm{ t\mapsto t^\al CT(t)x }_{L^p(\RR_+,Z)}
       \la^\al \bignorm{ (\la t)^{-\al} e^{-\la t} }_{L^{p'}(\RR_+)} \\
(s=\la t)
&\leq &  K  \la^{\al+\sfrac{1}{p'}} \; \norm{ x }_X 
\end{eqnarray*}
where  the number $K$ depends
only on $p$ and the norm of  $\Psi_\infty$. 
Next we treat the special case $p=\infty$.
\begin{proposition}\label{prop:p-infty}
Let $A$ be an injective sectorial operator of type $\om<\pihalbe$ 
on a Banach space $X$ and let $C \in B(X_1, Z)$. For $\al\in(0,1-\sfrac{1}p)$ the following
assertions are equivalent:
  \let\ALTLABELENUMI\labelenumi \renewcommand{\labelenumi}{(\roman{enumi})}  
  \let\ALTTHEENUMI\theenumi \renewcommand{\theenumi}{(\roman{enumi})}  
\begin{enumerate}
\item \label{item:p-infty-a} 
       The map $\Psi_\infty$ is bounded $X \to L^\infty_\al(\RR_+, Z)$
\item \label{item:p-infty-b}  
        The set $\{ \la^{1-\al} C (\la+A)^{-1}:\; \la>0\}$ is bounded in  $B(X, Z)$.
\item \label{item:p-infty-c} 
        $C$ is bounded in the norm $(X, \dot{X}_1)_{\al,1} \to Z$.
\end{enumerate}
  \let\labelenumi\ALTLABELENUMI    
  \let\theenumi\ALTTHEENUMI
\end{proposition}
\begin{proof}
From the necessary condition \eqref{eq:necess-C-zul} it follows directly that
\ref{item:p-infty-a} implies \ref{item:p-infty-b}.

Assume that \ref{item:p-infty-b} holds. For $x\in \DOMAIN(A)$ we have
\begin{equation}  \label{eq:triebel-gleichheit}
Cx = \int_0^\infty \la C A (\la+A)^{-2} x \,\tfrac{d\la}\la.
\end{equation}
Indeed, convergence follows from \ref{item:p-infty-b} since for small $\la>0$,
\[
\norm{\la^\al A(\la+A)^{-1}x} \leq c_1 \la^\al \norm{x}
\]
by sectoriality of $A$. For $\la\to\infty$, we have
\[
\norm{\la^\al A(\la+A)^{-1}x} \leq c_2 \la^{\al-1} \norm{A x} 
\]
also by sectoriality of $A$. Equality in \eqref{eq:triebel-gleichheit}
follows immediately from formula
\cite[1.14.2/(1)]{Triebel:interpolation}. Therefore, for $x \in
\DOMAIN(A)$, we obtain
\begin{eqnarray*}
     \norm{ Cx}_Z 
&\leq&c^{-1} \int_0^\infty \bignorm{\la^{1-\al} C (\la+A)^{-1}
            \la^{\al} A (\la+A)^{-1} x}_Z \,\tfrac{d\la}\la\\
&\leq&c^{-1} M \int_0^\infty \bignorm{\la^{\al} A (\la+A)^{-1} x}_X
           \,\tfrac{d\la}\la
\sim \norm{x}_{(X, \dot{X}_1)_{\al, 1}}, 
\end{eqnarray*}
where we used Proposition~\ref{prop:homogene-interpol-R}. So \ref{item:p-infty-b}
implies \ref{item:p-infty-c}.\\

Finally, let \ref{item:p-infty-c} hold. By
Proposition~\ref{prop:homogene-interpol-T} we then have
\begin{eqnarray*}
     \bignorm{t^\al C T(t) x}_Z 
&\leq&  M \, \bignorm{t^\al T(t)x}_{(X,\dot{X}_1)_{\al,1}} \\
&\sim& M \int_0^\infty \bignorm{t^\al s^{1-\al} A T(t+s)x}_X \,\tfrac{ds}s \\
&\leq&  \widetilde{M} \, \norm{x} \int_0^\infty t^\al s^{1-\al} (s{+}t)^{-1}\,\tfrac{ds}s \\
(s=\si t)
&=& 
       \widetilde{M} \, \norm{x} \int_0^\infty \si^{1-\al} (1+\si)^{-1} \, \tfrac{d\si}\si.
\end{eqnarray*}
The second estimate used the fact that, for bounded analytic semigroups,
the operators $(tA)\, T(t)$, $t>0$ are uniformly bounded.
\end{proof}
\begin{theorem}\label{thm:A1-p-q}
Let $1\leq p\leq q \leq \infty$ and $A$ be an injective sectorial 
operator of type $\om<\pihalbe$ on a Banach space $X$ and let 
$C \in B(X_1, Z)$. Let $\alpha\in(-\sfrac{1}p, 1-\sfrac{1}{p})$.  Then the following assertions hold:
\begin{enumerate}
\item \label{item:A1-p-q-a}
      If $\Psi_\infty$ is bounded $X \to L^p_\al(Z)$, then it
      is also bounded $X \to L^q_{\alpha+\sfrac{1}p-\sfrac{1}q}(Z)$.
\item \label{item:A1-p-q-b} 
      If $\Psi_\infty$ is bounded $L^q_{\alpha+\sfrac{1}p-\sfrac{1}q}(Z)$ and 
      if $X \emb (\dot{X}_{-1},\dot{X}_1)_{\einhalb,p}$,
      then it is also bounded $X \to L^p_\al(Z)$.
\end{enumerate}
\end{theorem}

Theorem~\ref{thm:discussion-A1-A3}~\ref{item:discussion-A1}
is a corollary of this result letting $q=\infty$. Before giving a
proof we point out its main argument, a simple reiteration observation.

\medskip
\begin{observation}\label{obs:observation}
Let numbers $p,q\in[1,\infty]$ and $\theta \in(0,1)$ be
given. Then,
\begin{equation}\label{eq:heart-of-the-matter}
  \begin{split}
\bigl(     \bigl( \dot{X}_{-1}, X \bigr)_{\theta,q}, 
           \bigl( X ,\dot{X}_{1}  \bigr)_{\theta,q} 
    \bigr)_{1-\theta,p} 
&=\bigl( 
           \bigl( \dot{X}_{-1}, \dot{X}_{1} \bigr)_{\sfrac{\theta}2,q}, 
           \bigl( \dot{X}_{-1}, \dot{X}_{1} \bigr)_{\einhalb+\sfrac{\theta}2,q}
    \bigr)_{1-\theta,p} \\
&= \bigl( \dot{X}_{-1}, \dot{X}_{1} \bigr)_{\einhalb,p}.
  \end{split}
\end{equation}
The first equality in (\ref{eq:heart-of-the-matter}) holds by
Lemma~\ref{lem:einbettung} and reiteration for the real method. The second
equality is the reiteration formula, see \cite[Theorem 1.10.2]{Triebel:interpolation}.
\end{observation}
\begin{proof}[{Proof of {Theorem~\ref{thm:A1-p-q}}}]
\ref{item:A1-p-q-a}. If $\Psi_\infty$ is bounded $X \to L^p_\al(Z)$,
then by the necessary condition \eqref{eq:necess-C-zul} and Proposition~\ref{prop:p-infty},
$\Psi_\infty$ is bounded  for $X \to L^\infty_{\al+\sfrac{1}p}(Z)$. Therefore,
\begin{eqnarray*}
&&     \biggl(\int_0^\infty \bignorm{ t^{\sfrac{1}p-\sfrac{1}q+\alpha}
       CT(t) x}^q \,dt\biggr)^{\sfrac{1}{q}} \\
&=&  \biggl(\int_0^\infty \bignorm{ t^\alpha CT(t) x}^p \; \bignorm{
       t^{\al+\sfrac{1}p} CT(t)x}^{q-p}\,dt 
       \biggr)^{\sfrac{1}{q}} \\
&\leq&\biggl(\int_0^\infty \bignorm{ t^\alpha CT(t) x}^p \,dt
       \biggr)^{\sfrac{1}{q}} \cdot
     \biggl( \sup_{t>0} \; t^{\al+\sfrac{1}p}
       \bignorm{CT(t)x} \biggr)^{1-\sfrac{p}{q}}\\
&\leq& c \, \bignorm{x}^{\sfrac{p}{q}} \cdot \bignorm{x}^{1-\sfrac{p}{q}}
= c \,\bignorm{x} 
\end{eqnarray*}
Now assume that \ref{item:A1-p-q-b} holds. By the necessary condition
\eqref{eq:necess-C-zul} and Proposition~\ref{prop:p-infty}, $C$ is
bounded in norm $E \to Z$, where $E := (X, \dot{X}_1)_{\theta,1}$ with
$\theta := \al{+}\sfrac{1}p$.  The part $A_E$ of $\AA_{-1}$ in $E$ is
injective and sectorial in $E$. Therefore, we can define $\dot{E}_{-1}
:= \dot{E}_1(A_E)$ and obtain $\dot{E}_{-1} =
(\dot{X}_{-1},X)_{\theta,1}$ whence $\bigl( \dot{E}_{-1}, E
\bigr)_{1-\theta,p} =\bigl( \dot{X}_{-1}, \dot{X}_{1}
\bigr)_{\einhalb,p}$ by letting $q=1$ in
Observation~\ref{obs:observation}. Taking this into account, we obtain
by $C\in B(E,Z)$ and Proposition~\ref{prop:homogene-interpol-T}
\begin{eqnarray*}
      \bignorm{ t\mapsto t^{\alpha} C T(t)x }_{L^p(\RR_+,Z)} 
&\leq& M \bignorm{ t\mapsto t^{\alpha} T(t)x }_{L^p(\RR_+,E)} \\
& = & M \bignorm{ t\mapsto t^{\al+\sfrac{1}{p}} T(t)x }_{L^p(\RR_+,dt/t,E)} \\
& = & M \bignorm{ t\mapsto t^{\al+\sfrac{1}{p}} A T(t)x }_{L^p(\RR_+,dt/t,\dot{E}_{-1})} \\
&\leq &\widetilde{M} \norm{x}_{(\dot{E}_{-1},E)_{1-\theta,p}} 
  =   \widetilde{M} \norm{x}_{(\dot{X}_{-1},\dot{X}_1)_{\einhalb,p}}
\end{eqnarray*}
Thus $\Psi_\infty$ is bounded $X \to L^p_\al(Z)$ since by assumption 
$X \emb (\dot{X}_{-1},\dot{X}_1)_{\einhalb,p}$.
\end{proof}

\subsection{Results on assumption \ref{item:B-zul}}
\label{sec:on-A2}

Theorem~\ref{thm:discussion-A1-A3}~\ref{item:discussion-A2} is in fact
covered by \cite[Theorems 1.8 and 1.9]{HaakKunstmann:weighted-Lp-admiss}. 
We 
shortly sketch the basic idea of the proof.
Necessity of the
boundedness of $B$ as stated in
Theorem~\ref{thm:discussion-A1-A3}~\ref{item:discussion-A2}
follows from
\begin{equation} \label{eq:necess-B-zul}
  \RANGE(B) \subseteq ( \dot X_{-1}, X )_{2(\al+\sfrac{1}p), \infty}  
.
\end{equation}
Indeed, consider the function $u(s) := \eins_{(\sfrac{t}2, t)}
u_0$. Then assumption \ref{item:B-zul} shows
\[
\bignorm{ A^{-\sfrac1p} [T(t)-T(\sfrac{t}2)] B u_0 }_{\dot X_{-1-\sfrac1p}} 
    \leq c_{\al, p} t^{\theta} \qquad t>0.
\]
with 
$\theta = 2\al +\sfrac1p$.
By \cite[Theorem 6.4.2]{Haase:Buch} 
we therefore also have
\[
\norm{ t^{-\theta} (tA)^{1-\sfrac1p} T(t) B u_0 }_{\dot X_{-1-\sfrac1p}} 
= 
\norm{ t^{1-\sfrac1p-\theta} A T(t) B u_0 }_{\dot X_{-1}}
\leq c_{\al, p}, \qquad t>0
\]
whence \eqref{eq:necess-B-zul} follows from
Proposition~\ref{prop:homogene-interpol-T}. Conversely, by analyticity of the
semigroup $T(\cdot)$, we have the pointwise estimates
\[
\norm{T_{-1}(t)}_{X\to X} \leq c
  \quad\text{and}\quad
\norm{T_{-1}(t)}_{\dot X_{-1} \to X} \leq c \,t^{-1}.
\]
Thus, by real interpolation, \eqref{eq:necess-B-zul} implies
$\norm{T(t)B}_{W \to X} \leq M t^{-\ga}$ for all $t>0$ where
$\ga = 1-2(\al +\sfrac1p)$ (see \cite[Lemma 1.12]{HaakKunstmann:weighted-Lp-admiss}).
Thus 
\begin{eqnarray*}
 \biggnorm{ \int_0^t T(t{-}s)Bu(s)\,ds }_X
 &\leq& \int_0^t \bignorm{ T(t{-}s)B }_{W\to X} \, s^{-\al} \, \norm{ s^\al u(s) }_W\,ds\\
 &\leq& c\int_0^t (t{-}s)^{-\ga} s^{-\al} \norm{ s^\al u(s) }_W\,ds.
\end{eqnarray*}
Let $k_{\al,\ga}(t,s)= \eins_{(0,t)}(s)(t{-}s)^{-\ga} s^{-\al}$ for $s,t\in(0,\tau)$.
Therefore, \ref{item:B-zul} is bounded if 
 the kernel $k_{\al,\ga}$ induces
a bounded integral  operator $K_{\al,\ga}:L^p(0,\tau)\to
L^\infty(0,\tau)$. This however follows from 
 following lemma that is taken from
\cite[Theorem 7]{HardyLittlewood:FractionalIntegrals}, see also
\cite[Theorem B]{SteinWeiss:FractionalIntegrals}).
\begin{lemma}[Hardy, Littlewood]\label{lem:HardyLittlewood}
Let $1< q < p \le\infty$ and $\ga\in(0,1)$. Then, for any numbers $0 \le \al < \beta < 1$ 
satisfying $1+\al - \beta -\ga = \sfrac1q -\sfrac1p >0$, the operator
\[
(T_\ga f)(t) := \int\limits_0^t \frac{f(s)}{(t-s)^\ga} \,ds
\]
is bounded  $L^q_\beta(\RR_+) \to L^p_\al(\RR_+)$. This also holds 
if $p=q=\infty$ or if $\al=\beta=0$.
\end{lemma}
The original proof of Hardy and Littlewood is incorrect (in 
\cite[displayed formula after (4.14) of p. 579]{HardyLittlewood:FractionalIntegrals}, 
see also a comment and a corrected proof in 
\cite[p. 504]{SteinWeiss:FractionalIntegrals}). We provide here a short 
interpolation argument.
\begin{proof}
(1) Let $p=q=\infty$. It suffices to verify $k(t,\cdot) \in L^1(\RR_+)$ with 
a uniform norm bound for $t\in (0,\tau)$ where 
$k(t,s) := \eins_{[0,t]}(s)(t{-}s)^{-\ga} t^\al s^{-\beta}$. A simple 
substitution shows that the characterising condition is 
\[
\beta, \ga < 1 \quad\text{and}\quad 1+\al = \beta+\ga.
\] 
(2) Next we consider the case $\al = \beta = 0$. 
Since 
$t^{-\ga} \in L^{\sfrac{1}\ga, \infty}$,
 a version of 
Young's inequality (see e.g. \cite[Theorem 1.2.13]{Grafakos})
yields $t^{-\ga}\ast: L^s \to L^r$ for 
$1+\sfrac{1}r = \sfrac{1}s + \ga$, $r,s \not\in \{1,\infty\}$.  \\
(3) The general case now follows by complex interpolation: 
(see e.g. \cite[Theorem 1.18.5]{Triebel:interpolation}) of (1) and (2):
\[
L^q_\beta = [L^s, L^\infty_{\sfrac\beta\theta}]_\theta 
  \quad\text{and}\quad
L^p_\al = [L^r, L^\infty_{\sfrac\al\theta}]_\theta 
\]
provided that $\sfrac1q = (1-\theta)\sfrac1s$ and 
$\sfrac1p = (1-\theta)\sfrac1r$. Moreover, by (2), 
\[
t^{-\ga}\ast: L^s \to L^r \qquad\text{holds for } 
   1+\sfrac{1}r = \sfrac{1}s + \ga\]
and by (1), 
\[
t^{-\ga}\ast: L^\infty_{\sfrac\beta\theta} \to L^\infty_{\sfrac\al\theta} 
  \qquad\text{holds for } \beta, \ga < \theta 
  \quad\text{and}\quad \theta+\al = \beta+\theta\ga.
\]
Under the assumptions on $\al,\beta, \ga, p,q$  all of the 
above conditions are satisfied.
\end{proof}

Thus, $K_{\al,\ga}$, given by 
$k_{\al,\ga}(t,s)= \eins_{(0,t)}(s)(t{-}s)^{-\ga} s^{-\al}$ for $s,t\in(0,\tau)$
is bounded $L^p(0,\tau)\to L^\infty(0,\tau)$ provided that
one of the following conditions holds:
\begin{equation}\label{eq:bedingungen-fuer-control-zul-typ-al}
\begin{array}{rlllll}
 (i)  & p=1 & \tau<\infty & \al\leq0 & \ga\leq0 & \\
 (ii) & p=1 & \tau=\infty & \al=0   & \ga=0   & \\
 (iii)& p>1 & \tau<\infty & \al+\sfrac{1}{p}<1 & \ga+\sfrac{1}{p}<1 
            & \al+\ga+\sfrac{1}{p}\leq1\\
 (iv) & p>1 & \tau=\infty & \al+\sfrac{1}{p}<1 & \ga+\sfrac{1}{p}<1 
            & \al+\ga+\sfrac{1}{p}=1.
\end{array}
\end{equation}
Observe that condition $(iv)$ implies $\al>0$ and $\ga>0$. 
In case $\al>0$ the assertion of 
Theorem~\ref{thm:discussion-A1-A3}~\ref{item:discussion-A2}
now follows immediately. In case $\al=0$ we follow a different strategy. 
Instead of analysing boundedness of the convolution operator $T(\cdot)\ast$ 
we can study boundedness of the map
\begin{equation} \label{eq:duale-zul-def}
  \Phi_\tau:  L^p_\al((0,\tau),W) \to X,
 \qquad
  \Phi_\tau(u) = \int_0^\tau T(t)B u(t)\,dt
\end{equation}
for $p\in [1, \infty]$ and $\al=0$. Therefore in some sense we are in 
a
dual situation to the discussion of assumption \ref{item:C-zul} and indeed
similar methods can be employed. We shall discuss boundedness of 
$\Phi_\tau$ in \eqref{eq:duale-zul-def} for general $\al$ and then 
deduce the remaining step for 
Theorem~\ref{thm:discussion-A1-A3}~\ref{item:discussion-A2} as a special case.
\begin{proposition}\label{prop:p-1}
Let $A$ be an injective sectorial operator of type $\om<\pihalbe$ on a 
Banach space $X$. Let $B \in B(W, X_{-1})$ and let $\al\in(0,1)$.
Then the following assertions are equivalent:
\let\ALTLABELENUMI\labelenumi \renewcommand{\labelenumi}{(\roman{enumi})}  
\let\ALTTHEENUMI\theenumi \renewcommand{\theenumi}{(\roman{enumi})}  
\begin{enumerate}
\item \label{item:p-1-a}
      The map $\Phi_\tau$ is bounded $L^1_\al(W) \to X$.
\item \label{item:p-1-b}
      The set $\{ \la^{1-\al} (\la+A_{-1})^{-1}B:\; \la>0\}$ is
      bounded in  $B(W,X)$.
\item \label{item:p-1-c}
      $B$ is bounded in the norm $W \to (\dot{X}_{-1},X)_{1-\al,\infty}$.
\end{enumerate}
  \let\labelenumi\ALTLABELENUMI
  \let\theenumi\ALTTHEENUMI
\end{proposition}
\begin{proof}
The implication \ref{item:p-1-a} $\Rightarrow$ \ref{item:p-1-b} is
similar to Proposition~\ref{prop:p-infty}.\\
Let \ref{item:p-1-b} hold, and let $u\in U$. By
Proposition~\ref{prop:homogene-interpol-R}, we have 
\begin{eqnarray*}
       \norm{Bu}_{(\dot{X}_{-1},X)_{1-\al,\infty}} 
&\sim& \sup_{\la>0} \bignorm{ \la^{1-\al} 
         A_{-1} (\la+A)^{-1} B u}_{\dot{X}_{-1}} \\
&=& \sup_{\la>0} \bignorm{ \la^{1-\al} (\la+A)^{-1} B u}_{X}
\leq K \, \norm{u}.
\end{eqnarray*}
Hence \ref{item:p-1-b} implies \ref{item:p-1-c}. 
Finally, if \ref{item:p-1-c} holds, then, by Proposition~\ref{prop:homogene-interpol-T},
\begin{eqnarray*}
   \biggnorm{ \int_0^\infty t^\al T_{-1}(t) B u(t) \,dt}_X
&\leq&\int_0^\infty \bignorm{t^\al A_{-1} T_{-1}(t)  B
    u(t)}_{\dot{X}_{-1}} \,dt \\
&\leq& \int_0^\infty \sup_{s>0} \bignorm{s^\al A_{-1} T_{-1}(s)B
    u(t)}_{\dot{X}_{-1}} \,dt \\
&\leq& c \int_0^\infty \bignorm{B
    u(t)}_{(\dot{X}_{-1},X)_{1-\al,\infty}} \,dt \\
& = & c \, \norm{u}_{L^1(\RR_+, U)} \, \norm{ B }_{U \to (\dot{X}_{-1},X)_{1-\al,\infty}},
\end{eqnarray*}
which proves the last implication.
\end{proof}

\begin{theorem}\label{thm:A2-p-q}
Let $1\leq r\leq p\leq \infty$ and $A$ be an injective sectorial 
operator of type $\om<\pihalbe$ on a Banach space $X$. Let $B \in
B(W,X_{-1})$ and let $\alpha\in(-\sfrac{1}{p'},
1-\sfrac{1}{p'})$. Then the following assertions hold:
\begin{enumerate}
\item \label{item:A2-p-q-a} 
      If $\Phi_\tau$ is bounded  $L^p_\al(W)\to X$, then it is also bounded
      $L^r_{\alpha+\sfrac{1}r-\sfrac{1}p}(W) \to X$.
\item \label{item:A2-p-q-b} 
      If $\Phi_\tau$ is bounded  $L^r_{\alpha+\sfrac{1}r-\sfrac{1}p}(W)
      \to X$  and if $(\dot{X}_{-1}, \dot{X}_{1})_{\einhalb,p} \emb X$,
      then it is also bounded $L^p_\al(W) \to X$.
\end{enumerate}
\end{theorem}
\begin{proof}
\ref{item:A2-p-q-a}. By assumption we have
\begin{equation*}
  \biggnorm{\int_0^\infty t^{\sfrac{1}r-\sfrac{1}p+\al} T_{-1}(t)B
    u(t) \,dt }
\leq c_1 \, \bignorm{ t^{\sfrac{1}r-\sfrac{1}p} u(t) }_{L^p(\RR_+,W)}.
\end{equation*}
Proposition~\ref{prop:p-1} shows that $\Phi_\tau$ is bounded $L^1_{\al+\sfrac{1}{p'}}(W) \to X$ whence
\begin{equation*}
  \biggnorm{\int_0^\infty t^{\sfrac{1}r-\sfrac{1}p+\al} T_{-1}(t)B
    u(t) \,dt }
\leq c_2 \, \bignorm{ t^{\sfrac{1}r-1} u(t) }_{L^1(\RR_+,W)}.
\end{equation*}
These two estimates allow interpolation by the complex method 
with $\theta \in (0,1)$ chosen such that  $\sfrac{1}r = \theta \cdot 1 +
(1-\theta)\sfrac{1}p$. Applying \cite[Theorem 5.5.3]{BerghLoefstroem},
one obtains
\begin{equation*}
  \biggnorm{\int_0^\infty t^{\sfrac{1}r-\sfrac{1}p+\al} T_{-1}(t)B
    u(t) \,dt }
\leq c_3 \, \bignorm{ t^{\tau} u(t) }_{L^r(\RR_+,W)},
\end{equation*}
where $\tau = \theta(\sfrac{1}r-1)+(1-\theta)(\sfrac{1}r-\sfrac{1}p) =
0$, and the assertion is proved.\\
\ref{item:A2-p-q-b}. If $\Phi_\tau$ is bounded
$L^r_{\alpha+\sfrac{1}r-\sfrac{1}p}(W) \to X$, then $B$ is bounded in
norm $W \to F$ where $F := (\dot{X}_{-1}, X)_{1-\si, \infty}$ with
$\si = (\al+\sfrac{1}{p'})/1$. 
Notice that $\dot{F}_{1} = (X,\dot{X}_{1})_{1-\si, \infty}$, and that we have 
$(F,\dot{F}_{1})_{\si,p} = (\dot{X}_{-1},\dot{X}_{1})_{\einhalb,p}$
by letting $q=\infty$ and $\theta=1-\si$ in Observation~\ref{obs:observation}.
By $(\dot{X}_{-1}, \dot{X}_{1})_{\einhalb,p} \emb X$ we thus have
\begin{eqnarray*}
&&      \biggnorm{ \int_0^\infty s^{\al} T_{-1}(s) B u(s)\,ds }_X \\
&\leq& c \, \biggnorm{ \int_0^\infty s^{\al}
  T_{-1}(s) B u(s)\,ds}_{(F,\dot{F}_{1})_{\si,p}} \\
&=&   c \, \biggnorm{t\mapsto \int_0^\infty s^{\al} t^{1-\si} A
  T_{-1}(s+t) B u(s)\,ds}_{L^{p}(\RR_+,dt/t,F)} \\
&=&   c \, \biggnorm{t\mapsto \int_0^\infty
      s^{\al} t^{1-\al-\sfrac{1}{p'}-\sfrac{1}p}
      A T_{-1}(s+t) B u(s)\,ds}_{L^{p}(\RR_+,F)}.
\end{eqnarray*}
Notice that the operator-valued kernel $K(s,t) :=
s^{\al} t^{-\al} AT_{-1}(s+t)$ satisfies 
\[
    \norm{K(s,t)} 
\leq M\, \frac{ s^{\al} t^{-\al}} {t+s} =: k(s,t)
\]
since $T(\cdot)$ is bounded analytic. The scalar
kernel $k(\cdot,\cdot)$ is homogeneous of degree
$-1$ and, by $\al \in (-\sfrac{1}{p'}, 1-\sfrac{1}{p'})$, 
the function 
\[
s \mapsto s^{-\sfrac{1}{p}}k(s,1) = \frac{s^{\al-\sfrac{1}p}}{1+s}
\]
is integrable over $\RR_+$. By \cite[Lemma A.3]{Stein:singular-int}, 
we thus obtain
\[
\biggnorm{ \int_0^\infty s^{\al} T_{-1}(s) B u(s)\,ds }_X 
\leq c'\norm{B u}_{L^p(\RR_+,F)}\leq c''\norm{u}_{L^p(\RR_+,W)},
\]
as desired.
\end{proof}

\section{Appendix}
\subsection*{Function spaces based on Lorentz spaces}\label{sec:lorentz-def-und-interpol}

It is well known (see e.g. \cite[1.18.6]{Triebel:interpolation}) that
real interpolation of Lebesgue spaces yields Lorentz spaces $L^{q,r}$
\begin{equation} \label{eq:lorentz-spaces}
\bigl( L^{q_0}, L^{q_1} \bigr)_{\theta, r} = L^{q,r}\qquad
\text{where } \tfrac{1}q = \tfrac{1-\theta}{q_0}+\tfrac{\theta}{q_1}
\end{equation}
and similar sequence spaces $\ell^{q, r} = (\ell^{q_0},
\ell^{q_1})_{\theta, r}$, see e.g. \cite[1.18.6 and 1.18.3]{Triebel:interpolation}.
In case $r=\infty$ these coincide with the weak Lebesgue (or Marcinkiewicz)
spaces $L^{q, \infty}$ for $q\in [1, \infty)$ with norm
\[
\norm{f}_{L^{q,\infty}(\Omega)} := \sup_{t >0} \; t \cdot \mu\bigl( \{ \omega \in \Omega:  \norm{f(\omega)} > t \}^\sfrac{1}q \bigr),
\]
and $\ell^{q,\infty}$ is the weak $\ell^q$ space of all sequences
$(c_\la)_{\la\in\NN}$ such that 
\[
   \sup_{\ga>0} \; \ga^q \,|\{\la\in\NN: |c_\la|>\ga\} |<\infty.
\]
In the sequel we study function spaces on $\RR^n$ constructed on
Lorentz spaces. When lifting (\ref{eq:lorentz-spaces}) with $(I-\Delta)^{-\sfrac{s}2}$ 
one obtains the spaces $H^s_{q,r}$ as corresponding real interpolation
spaces of $H^s_q$--spaces: 
\begin{equation}  \label{eq:weak-lifting}
\xymatrix{
   L^{q_0}   \ar@{--}[rr] \ar@{->}[d]_{(1-\Delta)^{-\sfrac{s}2}}  
&& L^{q,r} \ar@{--}[rr] \ar@{..>}[d]
&& L^{q_1} \ar@{->}[d]^{(1-\Delta)^{-\sfrac{s}2}}      \\
   H^s_{q_0} \ar@{--}[rr]
&& H^s_{q,r} \ar@{--}[rr]       
&& H^s_{q_1}
}
\end{equation}
We have thus  $\norm{x}_{H^s_{q,r}} \sim \norm{(1-\Delta)^{\sfrac{s}2}  x}_{L^{q,r}}$.

\smallskip

On the other hand, homogeneous and inhomogeneous Besov- and
Triebel-Lizorkin type spaces on $\RR^n$ may be 
constructed from Lorentz spaces by replacing the $L^q$--norm in the space variable by 
the corresponding Lorentz norm $\norm{\cdot}_{L^{q, r}}$. We denote these spaces by
$B^s_{(q,r), p}$ and $F^s_{(q,r), p}$ and 
in the homogeneous case by $\dot B^s_{(q,r), p}$ and
$\dot F^s_{(q,r), p}$
(in Triebel's book (see \cite[Definition
2.4.1]{Triebel:interpolation}) these spaces are denoted by $B^s_{q,p,  (r)}$ and 
$F^s_{q, p, (r)}$, respectively). 
It follows from a 'horizontal' interpolation property of
Triebel-Lizorkin spaces in $(\ast)$ below (see \cite[Theorem
2.4.2/5]{Triebel:interpolation}) 
that $H^s_{q, r} = F^s_{(q,r), 2}$ for $s \in\RR$, $p,q \in
(1,\infty)$ and $r\in [1,\infty]$. Indeed
\begin{eqnarray*}
F^s_{(q,r), 2} & \overset{(\ast)}{=} & ( F^s_{q_0, 2}, F^s_{q_1, 2} )_{\theta,r} 
             = ( H^s_{q_0}, H^s_{q_1} )_{\theta, r} \\
           & = & (I-\Delta)^{-\sfrac{s}2} ( L^{q_0}, L^{q_1} )_{\theta, r}
             = (I-\Delta)^{-\sfrac{s}2} L^{q,r}
             = H^s_{q, r}
\end{eqnarray*}
To define homogeneous spaces $\dot H^s_{q, r}$ we employ the same
technique as above but lift with $(-\Delta)^{-\sfrac{s}2}$ instead of
$(I-\Delta)^{-\sfrac{s}2}$. One obtains the analogue identity $\dot H^s_{q, r} = 
\dot F^s_{(q,r), 2}$. The proof is an immediate consequence
of a retraction/co-retraction argument that boils down the problem to
an interpolation of vector-valued Lebesgue spaces.

\medskip

For Besov type spaces a similar 'horizontal' interpolation property 
holds: let $s_0, s_1 \in \RR$,
$p_0, p_1 \in [1, \infty)$ and $q_0, q_1 \in (1, \infty)$ with
$q_0\not= q_1$. Then
\begin{equation}\label{eq:besov-horizontal}
\bigl( B^{s_0}_{q_0, p_0}, B^{s_1}_{q_1, p_1} \bigr)_{\theta, p} =
B^s_{(q,p), q}
\end{equation}
provided that $s=(1{-}\theta)s_0 + \theta s_1$ and $\tfrac1q =
\tfrac{1-\theta}{q_0}+\tfrac{\theta}{q_1}$
and $\tfrac1p = \tfrac{1-\theta}{p_0}+\tfrac{\theta}{p_1}$, see
\cite[Theorem 2.4.1/5]{Triebel:interpolation}. 
If 
$p\in [1,\infty]$ is not the interpolated index of $p_0$ and $p_1$, the above
equality may not hold. However, we always have the embedding property
(\ref{eq:einschachtelung-schwacher-raeume}) below. 
For its proof we require the following lemma.

\begin{lemma}\label{lem:claim-one}
Let $s\in\RR$, $q_0,q_1 \in (1,\infty)$. Then 
$B^s_{q_j,1}\emb H^s_{q_j}\emb B^s_{q_j,r}$, $j=0,1$
with $r=\max(q_0,q_1,2)$
\end{lemma}
\begin{proof}
For $q_j\geq 2$ this follows from
\[
H^s_{q_j} = F^s_{q_j,2} \emb F^s_{q_j, q_j} = B^s_{q_j,q_j} \emb B^s_{q_j, r}
\]
and for 
$q_j<2$
 this follows from 
$H^s_{q_j} = F^s_{q_j,2} \overset{(\ast)}{\emb} B^s_{q_j, 2} \emb
B^s_{q_j, r}$ where we use Minkowski's inequality in $(\ast)$. 
\end{proof}

Real interpolation of the embedding in Lemma~\ref{lem:claim-one} with $(\theta,
p)$ for fixed $s\in\RR$ and $q_0\not= q_1$ yields the following result 
due to {\sc Peetre} \cite[Theorem 1]{Peetre:weak-Besov} (see also 
\cite[Remark 4,Section 2.4.1]{Triebel:interpolation}) by (\ref{eq:besov-horizontal}). 
\begin{equation}
   \label{eq:einschachtelung-schwacher-raeume}
      B^s_{(q,p),\min(p,r^\ast)} 
  \emb H^s_{q,p} 
  \emb B^s_{(q,p),\max(p,r^\ast)}.
\end{equation}
Here, $r^*$ may chosen to satisfy 
$\tfrac{1}{r^\ast} = 1-\theta + \tfrac{\theta}r$ for any $r > \max(q, 2)$. 
Notice that 
\[
S: B^s_{(q,r),p} \to \ell_p^s(L^{q,r})  \qquad  S f := ( f \ast \phi_j),
\quad j\in \NN
\]
\[
\dot S: \dot B^s_{(q,r),p} \to \ell_p^s(L^{q,r})  \qquad  \dot S f := ( f
\ast \phi_j) \quad j\in \ZZ
\]
are co-retractions in the sense of 
\cite[Definition 1.2.4]{Triebel:interpolation} (see
also \cite[(2.3.2/12)]{Triebel:interpolation} for more
details). This means that (\ref{eq:besov-horizontal}) and also
the above proof 
transfer to the case of homogeneous spaces
as well, since they rely essentially on an embedding result for
interpolation spaces of vector-valued $\ell_q^s$--spaces. 
The authors thank {\sc H. Triebel} for suggesting the co-retraction argument.
\begin{lemma}\label{lem:claim-two}
One has $\dot H^{-1}_{\sfrac{q}2,\infty} \emb \dot B^t_{(q,\infty),\infty}$ provided that $-1-\sfrac{n}q = t$.
\end{lemma}
\begin{proof} Notice that by the well-known Sobolev embedding,  $H^s_q
\emb L^{q^\ast}$ where $s > 0$ and $\tfrac{n}{q^\ast}=\tfrac{n}q - s$. 
A scaling argument (i.e. regarding norm estimates for $u(\la \cdot)$
with $\la>0$) yields the estimate
\[
    \la^{-\sfrac{n}{q^\ast}} \norm{u}_{q^\ast} 
\leq C (\la^{-\sfrac{n}q} \norm{u}_q + \la^{s-\sfrac{n}q} \norm{\cF^{-1}(\xi \mapsto |\xi|^s\cF(u) }_q)
\]
Now multiplication with $\la^{\sfrac{n}{q^\ast}}=\la^{\sfrac{n}{q}-s}$ and letting $\la \to
\infty$ gives the embedding of the homogeneous space $\dot H^s_q \emb L^{q^\ast}$.
In particular, we have $\dot H^{\sfrac{n}{q_j}}_{\sfrac{q_j}2} \emb
L^{q_j}$. Real interpolation of this embedding for adequate values
$q_0$, $q_1$ and $\theta$ yields
\begin{equation}\label{eq:homogene-einbettung}
\dot H^{\sfrac{n}q}_{(\sfrac{q}2, \infty)} 
= (\dot H^{\sfrac{n}q_0}_{\sfrac{q_0}2}, \dot  H^{\sfrac{n}q_1}_{\sfrac{q_1}2} )_{\theta, \infty} 
\emb (L^{q_0}, L^{q_1})_{\theta,\infty} = L^{q,\infty}
\end{equation}
Now apply \eqref{eq:einschachtelung-schwacher-raeume} with $s=0$ in
$(\ast)$ to conclude
$\dot H^{\sfrac{n}q}_{(\sfrac{q}2, \infty)} \emb L^{q,\infty} \overset{(\ast)}{\emb} \dot B^0_{(q,\infty),\infty}$. 
Finally, lifting by $-\sfrac{n}q-1$ proves the lemma.
\end{proof}

\medskip

\textbf{Wavelet characterisation of Besov spaces.}
The homogeneous Besov space $\dot{B}^s_{q,p}(\RR^n)$ can be characterised
as the space of all wavelet series $\sum\al(\la)\psi_\la(x)$ such that 
\begin{equation}\label{eq:besov-meyer-char}
      \bigg(2^{js}\,2^{\sfrac{-nj}{q}}
           \bigg(\sum_{\la\in\Lambda_j}|\al(\la)|^q\bigg)^{\sfrac1q}
      \bigg)_j
      \in \ell^p(\ZZ),
\end{equation}
see \cite[page 198]{Meyer:ondelettes} if we normalise the $L^1$-norm of the
functions $\psi_\la$ (see \cite[Lemma 4.2.5]{Cannone:Buch}). By real interpolation
$(\cdot,\cdot)_{\theta,\infty}$ between $\dot{B}^s_{q_0,p}$ and
$\dot{B}^s_{q_1,p}$ we thus obtain a characterisation of
$\dot{B}^s_{(q,\infty),p}(\RR^n)$ in terms of wavelet coefficients
which reads as \eqref{eq:besov-meyer-char} but with
$\norm{(\al(\la))_{\la\in\Lambda_j}}_{\ell^{q,\infty}}$ replacing the
$\ell^q$-norm. 
We are interested in the case $s=-1-\sfrac{n}q$ and characterise
distributions in $\dot{B}^{-1+\sfrac{n}q}_{(q,\infty),\infty}(\RR^n)$ that are
homogeneous of degree $-1$. Observe that the characterisation in this case
reads
\begin{equation}\label{eq:spec-weak-besov-char}
     \norm{(\al(\la))_{\la\in\Lambda_0}}_{\ell^{q,\infty}}<\infty
\end{equation} 
since $js-\sfrac{nj}q=-j+\sfrac{nj}q-\sfrac{nj}q=-j$ and, by homogeneity of degree $-1$, 
the $\ell^{q,\infty}$-norm of $(\al(\la))_{\la\in\Lambda_j}$ equals 
$2^j\norm{(\al(\la))_{\la\in\Lambda_0}}_{\ell^{q,\infty}}$
(see \cite[(4.20)]{Cannone:Buch}). Repeating the arguments on
\cite[page 154/5]{Cannone:Buch} we thus see that a distribution that is
homogeneous of degree $-1$ belongs to
$\dot{B}^{-1+\sfrac{n}q}_{(q,\infty),\infty}$ if and only if its restriction
to the annulus $\Om:=\{x\in\RR^n:\einhalb<|x|< \sfrac32\}$ belongs to the inhomogeneous
Besov-like space $B^{-1+\sfrac{n}q}_{(q,\infty),(q,\infty)}(\Om)$ which can be
obtained, e.g. by real interpolation $(\cdot,\cdot)_{\theta,\infty}$
between the Besov spaces $B^{s}_{q_0,q_0}(\Om)$ and $B^s_{q_1,q_1}(\Om)$. 
Arguing as in the proof of \cite[Theorem 4.2.2]{Cannone:Buch} another
equivalent condition is that the restriction to the unit sphere
$S^{n-1}$ in $\RR^n$ belongs to
$B^{-1+\sfrac{n}q}_{(q,\infty),(q,\infty)}(S^{n-1})$. This allows to
construct examples by homogeneous extension. We remark that a
Littlewood-Paley characterisation of
$B^s_{(q,\infty),(q,\infty)}(\RR^{n-1})$ would read $f\in
B^s_{(q,\infty),(q,\infty)}(\RR^{n-1})$ if and only if
\[
  \norm{S_0f}_{L^{q,\infty}}
+\bigg\|( 2^{js}\|\Delta_jf\|_{L^{q,\infty}})_{j\ge0}\bigg\|_{\ell^{q,\infty}}
<\infty.
\] 

\subsection*{Morrey spaces and Sobolev-type spaces based on Morrey spaces}
Let $q\in(1,\infty)$ and $\la\in [0,\sfrac{n}q]$. Then the Morrey space
$\Morrey^{q,\la}(\RR^n)$ consists of all functions 
$f \in L^q_\text{loc}(\RR^n)$ for which the maximal function 
\[
M_{q,\la} f: \quad x \mapsto \inf_{c\in\RR} \; \sup_{r>0} \quad r^{\la} \biggl(
\tfrac{1}{|B(x, r)|} \, \int_{B(x,r)} \bigl| f(y) - c|^q \,dy \biggr)^\sfrac1q
\]
is bounded on $\RR^n$. Notice that the exponent of the radius $r$ associated to
the value of $\la$ in the notation $\Morrey^{q,\la}(\RR^n)$ is not
consistent for all publications on the subject, see e.g. 
\cite{Campanato:definition-Morrey-Camp-Raum}. The above notation
seems the most natural to us. In {\sc Peetre} \cite{Peetre:Morrey} 
the notation ${\mathscr E}^{-\la, q}$ is used for $\Morrey^{q,\la}(\RR^n)$

We briefly summarise some results on Morrey spaces:
it is clear that 
$\Morrey^{q,\sfrac{n}q}(\RR^n) = L^q(\RR^n)$. 
A celebrated result of {\sc John} and {\sc Nirenberg} \cite{JohnNirenberg} states
that 
$\Morrey^{q,0}(\RR^n) = \text{BMO}(\RR^n)$. 
Recall that in case 
$\la\in (0,\sfrac{n}q]$ which interests us, one may let $c=0$ in the
above definition. In view of interpolation techniques one has to
remark that the spaces $\Morrey^{q,\la}(\RR^n)$ do not interpolate
'horizontally' for $n\ge 2$, i.e. one has
\[
\bigl(\Morrey^{q_0, \la_0}, \Morrey^{q_1, \la_1}
\bigr)_{\theta, q} \subsetneq \Morrey^{q, \la}.
\]
when $\tfrac1q = \tfrac{1{-}\theta}{q_0} + \tfrac\theta{q_1}$ and 
$\tfrac1\la = \tfrac{1{-}\theta}{\la_0} + \tfrac\theta{\la_1}$,
see {\sc Blasco}, {\sc Ruiz} and {\sc Vega}
\cite{BlascoRuizVega:Morrey} for details and  more complete
references. For $s\in\RR$ we now denote
\begin{eqnarray*}
\dot\Morrey^{q,\la,s}(\RR^n) 
&=&\bigl \{f \suchthat  \sF^{-1}(\xi\mapsto|\xi|^s \sF f(\xi))
          \in \Morrey^{q,\la}(\RR^n) \bigr\}   \quad\text{and} \\
 \Morrey^{q,\la,s}(\RR^n) 
&=& \bigl\{f \suchthat  \sF^{-1}(\xi\mapsto(1+|\xi|^2)^{s/2} \sF f(\xi))
          \in \Morrey^{q,\la}(\RR^n) \bigr\}
\end{eqnarray*}
with the same range of $\la \in (0, \sfrac{n}q]$. These spaces are
homogeneous and inhomogeneous Sobolev-type spaces based on Morrey spaces. 
We recall the fact that the usual Calder\'on-Zygmund operators are
bounded on Morrey spaces, see {\sc Peetre}
\cite{Peetre:convolution-op-in-Morrey-spaces}, see also
\cite{ChiarenzaFrasca}. In particular, the Mihlin 
(sometimes also transcribed as Mikhlin or Michlin)
theorem on Fourier multipliers holds in Morrey spaces. Therefore, 
we have for $s\in\NN$:
\[
 \Morrey^{q,\la,s}(\RR^n)=\bigl\{f\in\Morrey^{q,\la}(\RR^n) \suchthat \forall \; |\al|\le
 s: \quad \partial^\al f\in\Morrey^{q,\la}(\RR^n)\,\bigr\}.
\]
Moreover, considering the heat semigroup $T(\cdot)$ in $\Morrey^{q,\la}(\RR^n)$, we
have by translation invariance of the space that $T(\cdot)$ is bounded
analytic in $\Morrey^{q,\la}(\RR^n)$. Denoting the generator by $\Delta$, we
then have $D(\Delta)=\Morrey^{q,\la,2}(\RR^n)$ and, still by
Calder\'on-Zygmund theory, 
$(\Morrey^{q,\la})^{\displaystyle{\cdot}}_{s,-\Delta}(\RR^n)
=\dot\Morrey^{q,\la,2s}(\RR^n)$ for any $s\in\RR$.

{\em Riesz potential operator and embeddings.} 
Let $I_s$ be the Riesz potential operator given by the convolution
kernel $|x|^{s-n}$. Then the following Sobolev inequality for Morrey
spaces holds. 

\begin{proposition}[{\cite[Theorem 3.1]{Adams:Riesz-potentials-on-Morrey-spaces}}] 
\label{prop:sobolev-emb-morrey}
Let $s>0$, $\nu \in (0, \sfrac{n}p]$, and  $s \in (0, \nu)$ such that
$\sfrac1r = \sfrac1p - \tfrac{s}{\nu p}$.
Then $I_s: \Morrey^{p,\nu} \to \Morrey^{r,\frac{\nu p}r}$  is a bounded
linear operator, or, equivalently,
\[
\dot{\Morrey}^{p,\nu,s} \emb {\Morrey}^{r,\frac{\nu p}r}.
\]
\end{proposition}

This is used in Subsection \ref{sec:morrey}.

\providecommand{\bysame}{\leavevmode\hbox to3em{\hrulefill}\thinspace}

\end{document}